\documentclass[a4paper,openany,11pt]{memoir}

\usepackage{amsmath, amsthm, amssymb, amsfonts}
\usepackage[labelfont={sf}]{caption}
\usepackage{cite}
\usepackage{colonequals}
\usepackage[autostyle]{csquotes}
\usepackage{enumerate}
\usepackage{enumitem}
\usepackage{euscript}
\usepackage{fancyhdr}
\usepackage{float}
\usepackage[T1]{fontenc}
\usepackage{framed}
\usepackage[utf8]{inputenc}
\usepackage[letterpaper,hmargin=1in,vmargin={0.75in,0.2in}]{geometry}
\usepackage{lmodern}
\usepackage{hyperref}
\usepackage{mathrsfs}
\usepackage{mathtools}
\usepackage{multicol}
\usepackage{parskip}
\usepackage{relsize}
\usepackage{stmaryrd}
\usepackage[svgnames]{xcolor}
\usepackage{thmtools}
\usepackage{tikz}
\usepackage{titlesec}
\usepackage{xfrac}
\usepackage{yfonts}

\newcommand{\arr}{\mathbb{R}}
\newcommand{\cee}{\mathbb{C}}

\newcommand{\deee}{\,\mathrm{d}}

\newcommand{\enn}{\mathbb{N}}
\newcommand{\ennplus}{\enn_+}
\newcommand{\rzee}{{\mathbb{R}/\mathbb{Z}}}

\newcommand{\zee}{\mathbb{Z}}

\DeclareMathOperator{\ach}{\mathbf{H}}
\DeclareMathOperator{\achn}{\mathbf{H}_\ast}

\DeclareMathOperator{\einf}{\mathbf{EInf}}

\DeclareMathOperator{\aee}{\mathrm{ae}}
\DeclareMathOperator{\bbr}{\mathrm{BBor}}

\DeclareMathOperator{\pol}{\mathrm{Pol}}
\DeclareMathOperator{\prob}{\mathrm{Prob}}

\DeclareMathOperator{\re}{\mathrm{Re}}

\DeclareMathOperator{\ret}{\mathtt{Ret}}
\DeclareMathOperator{\tor}{\mathtt{Tor}}
\DeclareMathOperator{\eir}{\mathtt{Eir}}

\DeclareMathOperator{\spn}{\mathrm{spn}}

\DeclareMathOperator{\verr}{\hspace{-0.2 pt} ; \hspace{-0.2 pt}}

\renewcommand{\dotfill}{%
  \leavevmode
  \cleaders\hbox to 0.65em{\hss\textbf{.}\hss}\hfill\kern0pt
}
\renewcommand{\epsilon}{\varepsilon}
\renewcommand{\hat}{\widehat}
\renewcommand{\phi}{\varphi}

\newcommand{\ddotfill}{%
  \leavevmode
  \cleaders\hbox to 0.65em{\hss.\hss}\hfill\kern0pt
}

\newcommand{\lvertt}{\left \vert}
\newcommand{\rvertt}{\right \vert}

\newcommand{\eusc}[1]{\EuScript{#1}}

\newcommand{\ov}[1]{\overline{#1}}

\makeatletter \def\vhrulefill#1{\leavevmode\leaders\hrule\@height#1\hfill \kern\z@}
\makeatother

\makeatletter
\newcommand{\hooktwoheadright}[2][]{%
  \lhook\joinrel
  \ext@arrow 0359\rightarrowfill@ {#1}{#2}%
  \mathrel{\mspace{-15mu}}\rightarrow
}
\makeatother

\newsavebox{\depbox}
\newlength{\deplabelwidth}
\newlength{\depboxwidth}
\titleformat{\part}[display]
  {\normalfont\sffamily\HUGE\bfseries}
  {Part \thepart}{20pt}{\HUGE}

\titleformat{\chapter}[display]
  {\normalfont\sffamily\LARGE\bfseries}
  {\chaptertitlename\ \thechapter}{20pt}{\Huge}
  
\titleformat{\section}
  {\normalfont\sffamily\LARGE\bfseries}
  {\thesection}{1em}{}
  
\titleformat{\subsection}
  {\normalfont\sffamily\Large\bfseries}
  {\thesubsection}{1em}{}
  
\titleformat{\subsubsection}
  {\normalfont\sffamily\Large\bfseries}
  {\thesubsubsection}{1em}{}

\definecolor{gloryred}{HTML}{BB133E}
\definecolor{gloryblue}{HTML}{2C3ABC}
\definecolor{glorygreen}{HTML}{24802A}
\definecolor{gloryorange}{HTML}{FF6C0C}
\definecolor{topgrey}{HTML}{CFCFCF}
\definecolor{midgrey}{HTML}{AAAAAA}
\definecolor{downgrey}{HTML}{828282}

\colorlet{headcolor1}{Black}
\colorlet{rulecolor1.0}{gloryred}
\colorlet{rulecolor1.1}{gloryblue}
\colorlet{rulecolor1.2}{glorygreen}

\newenvironment{leftbar1.0}{
  
  \MakeFramed {\advance\hsize-\width \FrameRestore}}
 {\endMakeFramed}
 
 \newenvironment{leftbar1.1}{
  
  \MakeFramed {\advance\hsize-\width \FrameRestore}}
 {\endMakeFramed}
 
 \newenvironment{leftbar1.2}{
  
  \MakeFramed {\advance\hsize-\width \FrameRestore}}
 {\endMakeFramed}

 \declaretheoremstyle[headfont=\sffamily\bfseries,
 notefont=\sffamily\bfseries,
 notebraces={}{},
 headpunct=,
 bodyfont=\sffamily,
 headformat=\color{headcolor1}\NAME~\NUMBER:~\!\!\NOTE\hfill\smallskip\linebreak,
 preheadhook=\begin{leftbar1.0},
 postfoothook=\end{leftbar1.0},
 ]{customTheorem}
 \declaretheorem[style=customTheorem, numberwithin=chapter]{Theorem}

 \declaretheoremstyle[headfont=\sffamily\bfseries,
 notefont=\sffamily\bfseries,
 notebraces={}{},
 headpunct=,
 bodyfont=\sffamily,
 headformat=\color{headcolor1}\NAME~\NUMBER:~\!\!\NOTE\hfill\smallskip\linebreak,
 preheadhook=\begin{leftbar1.1},
 postfoothook=\end{leftbar1.1},
 ]{customLemma}
 \declaretheoremstyle[headfont=\sffamily\bfseries,
 notefont=\sffamily\bfseries,
 notebraces={}{},
 headpunct=,
 bodyfont=\sffamily,
 headformat=\color{headcolor1}\NAME~\NUMBER\hfill\smallskip\linebreak,
 preheadhook=\begin{leftbar1.2},
 postfoothook=\end{leftbar1.2},
 ]{customProp}

 \declaretheoremstyle[headfont=\sffamily\bfseries,
 notefont=\sffamily\bfseries,
 notebraces={}{},
 headpunct=,
 bodyfont=\sffamily,
 headformat=\color{headcolor1}Proposition~\NUMBER:~\!\!\NOTE\hfill\smallskip\linebreak,
 preheadhook=\begin{leftbar1.2},
 postfoothook=\end{leftbar1.2},
 ]{customPropn}
 \declaretheorem[style=customPropn, numberwithin=section]{Propositionn}

\declaretheorem[style=customLemma, sibling=Propositionn]{Lemma}

\colorlet{headcolor2}{Black}
\colorlet{rulecolor2.0}{topgrey}
\colorlet{rulecolor2.1}{midgrey}
\colorlet{rulecolor2.2}{downgrey}

\newenvironment{leftbar2.0}{%
  \MakeFramed {\advance\hsize-\width \FrameRestore}}%
 {\endMakeFramed}
 
 \newenvironment{leftbar2.1}{%
  \MakeFramed {\advance\hsize-\width \FrameRestore}}%
 {\endMakeFramed}
 
 \newenvironment{leftbar2.2}{%
  \MakeFramed {\advance\hsize-\width \FrameRestore}}%
 {\endMakeFramed}

  \declaretheoremstyle[headfont=\sffamily\bfseries,%
 notefont=\sffamily\bfseries,%
 notebraces={}{},%
 headpunct=,%
 bodyfont=\sffamily,%
 headformat=\color{headcolor2}\NAME~\NUMBER:~\!\!\NOTE{}\hfill\smallskip\linebreak,%
 preheadhook=\begin{leftbar2.0},%
 postfoothook=\end{leftbar2.0},%
 ]{customDefinition}
 
\declaretheorem[style=customDefinition, sibling=Propositionn]{Definition}

 \declaretheoremstyle[headfont=\sffamily\bfseries,%
 notefont=\sffamily\bfseries,%
 notebraces={}{},%
 headpunct=,%
 bodyfont=\sffamily,%
 headformat=\color{headcolor2}\NAME~\NUMBER:~\!\!\NOTE{}\hfill\smallskip\linebreak,%
 preheadhook=\begin{leftbar2.1},%
 postfoothook=\end{leftbar2.1},%
 ]{customCorollary}

 \declaretheoremstyle[headfont=\sffamily\bfseries,%
 notefont=\sffamily\bfseries,%
 notebraces={}{},%
 headpunct=,%
 bodyfont=\sffamily,%
 headformat=\color{headcolor2}\NAME~\NUMBER:~\!\!\NOTE{}\hfill\smallskip\linebreak,%
 preheadhook=\begin{leftbar2.2},%
 postfoothook=\end{leftbar2.2},%
 ]{customRemark}

\declaretheorem[style=customRemark, sibling=Propositionn]{Remark}

\colorlet{headcolor3}{Black}
\colorlet{rulecolor3}{Black}

\newenvironment{leftbar3}{%
  \MakeFramed {\advance\hsize-\width \FrameRestore}}%
 {\endMakeFramed}

 \declaretheoremstyle[headfont=\sffamily\bfseries,%
 notefont=\sffamily\bfseries,%
 notebraces={}{},%
 headpunct=,%
 bodyfont=\sffamily,%
 headformat=\color{headcolor3}\NAME~\!\!\NOTE\hfill\smallskip\linebreak,%
 preheadhook=\begin{leftbar3},%
 postfoothook=\end{leftbar3},%
 qed=$\blacksquare$,
 ]{prf1}

\theoremstyle{Proof}
 
\declaretheorem[style=prf1, numberwithin=chapter]{Proof}

\colorlet{headcolor4}{Black}
\colorlet{rulecolor4}{gloryorange}

\newenvironment{leftbar4}{
  
  \MakeFramed {\advance\hsize-\width \FrameRestore}}
 {\endMakeFramed}

 \declaretheoremstyle[headfont=\sffamily\bfseries,
 notefont=\sffamily\bfseries,
 notebraces={}{},
 headpunct=,
 bodyfont=\sffamily,
 headformat=\color{headcolor4}Theorem~\NUMBER:~\!\!\NOTE\hfill\smallskip\linebreak,
 preheadhook=\begin{leftbar4},
 postfoothook=\end{leftbar4},
 ]{customoTheorem}
 \declaretheorem[style=customoTheorem, sibling=Propositionn]{oTheorem}
 
 \declaretheoremstyle[headfont=\sffamily\bfseries, notefont=\sffamily\bfseries, notebraces={}{}, headpunct=, bodyfont=\sffamily,
 headformat=\color{headcolor4}Theorem~0.1:~\!\!\NOTE\hfill\smallskip\linebreak,
 preheadhook=\begin{leftbar4}, postfoothook=\end{leftbar4}]{customoTheoremMaster}

\makeatletter
\newcommand{\forcedlabel}[2]{\phantomsection\def\@currentlabel{#1}\label{#2}}
\makeatother
 
 \declaretheoremstyle[headfont=\sffamily\bfseries,
 notefont=\sffamily\bfseries,
 notebraces={}{},
 headpunct=,
 bodyfont=\sffamily,
 headformat=\color{headcolor4}Lemma~\NUMBER:~\!\!\NOTE\hfill\smallskip\linebreak,
 preheadhook=\begin{leftbar4},
 postfoothook=\end{leftbar4},
 ]{customoLemma}
 \declaretheorem[style=customoLemma, sibling=Propositionn]{oLemma}

 \declaretheoremstyle[headfont=\sffamily\bfseries,
 notefont=\sffamily\bfseries,
 notebraces={}{},
 headpunct=,
 bodyfont=\sffamily,
 headformat=\color{headcolor4}Proposition~\NUMBER:~\!\!\NOTE\hfill\smallskip\linebreak,
 preheadhook=\begin{leftbar4},
 postfoothook=\end{leftbar4},
 ]{customoProp}

\numberwithin{equation}{section}

\fancypagestyle{plain}{%
\fancyhf{} 
\fancyfoot[C]{\sffamily \bfseries \fontsize{11pt}{11pt}\selectfont\thepage} 

}
\DeclareMathAlphabet{\mathpzc}{OT1}{pzc}{m}{it}

\DeclareQuoteStyle{single}{\lq}[]{\rq}{}[]{}
\DeclareQuoteStyle{double}{\lq\lq}[]{\rq\rq}{}[]{}

\thanksmarkseries{arabic}

\begin{document}

\sffamily

\setquotestyle{single}

\title{\Huge \bfseries Entropic rates from eigenfunction structure for $\times 2 \times 3$}

\author{\LARGE Peter Burton\thanks{\url{peterburton1728@gmail.com}} and Kate Juschenko\thanks{\url{k.juschenko@austin.utexas.edu}}}
\date{\LARGE\today}

\maketitle

\renewcommand{\abstractnamefont}{\Large \sffamily\bfseries}
\renewcommand{\abstracttextfont}{\sffamily}

\begin{abstract} In a companion article \cite{RJnodyn} the present authors proved a form of the Rudolph--Johnson theorem which is effective at finite resolution, bounding the low-frequency Fourier coefficients of a $\times 2 \times 3$-invariant measure by its normalized entropy. Here we read that bound in the other direction, as a constraint on eigenfunction behavior in potential counterexamples to Furstenberg's $\times 2 \times 3$ conjecture. We show that a near return of the coordinate function $x \mapsto e^{2\pi ix}$ to itself under the semigroup generated by $2$ and $3$ forces the resolution entropy to be small, and that adherence of the coordinate function to eigenfunctions of $x \mapsto 2x$ forces it to decay. These rates are explicit in the rate of adherence and in the arithmetic of the eigenvalues. \end{abstract}

\renewcommand\printtoctitle[1]{\LARGE\bfseries #1}

\renewcommand*{\cftsectionfont}{\sffamily}
\renewcommand*{\cftsectionpagefont}{\sffamily}

\renewcommand*{\cftsubsectionfont}{\sffamily}
\renewcommand*{\cftsubsectionpagefont}{\sffamily}

\renewcommand*{\tocheadstart}{\vspace*{12pt}}
\renewcommand*{\aftertoctitle}{\par\nobreak\vskip 6pt}

\tableofcontents*

\newpage

\section{Invitation to article}

This article is an investigation of potential counterexamples to Furstenberg's $\times 2 \times 3$ conjecture, the statement of which we now recall. Letting $\lambda$ denote Lebesgue measure on $\rzee$, Furstenberg asked in \cite{MR0213508} whether $\lambda$ is the unique atomless Borel probability measure on $\rzee$ which is invariant under both $x \mapsto 2x$ and $x \mapsto 3x$. Probably the key progress so far on this conjecture is the work of Rudolph \cite{MR1062766} and Johnson \cite{MR1194793}, which shows that $\lambda$ is in fact the unique ergodic invariant measure of positive entropy.\\

The present article is a companion to \cite{RJnodyn} by the same authors which proves a form of the Rudolph--Johnson theorem that is effective at finite resolution, bounding the low-frequency Fourier coefficients of a $\times 2 \times 3$-invariant measure in terms of its Shannon entropy on a partition of $\rzee$ into equal arcs. That bound is restated as Lemma \ref{cor.master} below.\\

Looking at the theory in the other direction, in this article  we show that a large Fourier coefficient at a low frequency forces the resolution entropy to be small. We exhibit two related spectral sources of large Fourier coefficients: near returns of the coordinate function $x \mapsto e^{2\pi i x}$ under the semigroup generated by $2$ and $3$, and adherence of the coordinate function to eigenfunctions of $x \mapsto 2x$.\\

The theorems that result describe how fast the resolution entropy of a $\times 2 \times 3$-invariant measure must decay once its Kronecker part carries the coordinate function at a quantified rate. If the answer to Furstenberg's question is affirmative, as is conjectured, then the measures satisfying these hypotheses are the atomic ones and the conclusions hold for trivial reasons. However, our theorems should be interpreted as constraints on the class of potential counterexamples to Furstenberg's conjecture and their content lies in the rates, which are explicit in every parameter.\\

We briefly introduce some notation needed to state four illustrative theorems.

\begin{itemize}
\item For $\theta \in \arr$ we adopt the standard notation $e(\theta) \coloneqq e^{2\pi i \theta}$. We identify $e$ with its descent to a function on $\rzee$, and refer to this as the \emph{coordinate function}.  
\item For $n \in \ennplus$ we consider the endomorphism $x \mapsto nx$ of $\rzee$. We define a Borel probability measure $\mu$ on $\rzee$ to be $\times n$\emph{-invariant} if the following holds for all continuous functions $\psi:\rzee \to \cee$.
\begin{align} \label{prose.intro.1} \int_{\rzee} \psi(nx) \deee \mu(x) = \int_{\rzee}  \psi(x) \deee \mu(x)  \end{align}
\item We identify $\rzee$ with its coset transversal $[0,1)$, and for $r \in \ennplus$ we make the following associated conventions. \begin{itemize} \item[$\circ$] We write $\ach(\mu,r)$ for the \emph{Shannon entropy} of the distribution induced by $\mu$ on the partition of $\rzee$ into $r$ arcs of equal length starting from $0$. 
\item[$\circ$] For $r \in \enn$ with $r \geq 2$ we adopt the notation $\achn(\mu,r) \coloneqq \ach(\mu,r)/\log(r)$ for the \emph{normalized Shannon entropy}, and note that this quantity always lies in $[0,1]$. 
\item[$\circ$] We refer the reader to Definition \ref{def.entropy} below for entropy-related formalisms.\end{itemize}
\item Given a $\times 2$-invariant Borel probability measure $\mu$ on $\rzee$, we define a nonzero $\psi \in L^2(\rzee,\mu)$ to be a \emph{$\times 2$-eigenfunction} with \emph{eigenvalue} $\zeta \in \cee$ if $\psi(2x) = \zeta \psi(x)$ for $\mu$-almost every $x \in \rzee$, and we define such $\psi$ to be \emph{normalized} if $\psi$ is a unit vector in $L^2(\rzee,\mu)$. We note that the invariance of $\mu$ guarantees $|\zeta|=1$.
\item For $x \in (1,\infty)_{\arr}$ we write $\log^{(1)}(x) \coloneqq \log(x)$, and for $m \in [2,\infty)_{\zee}$ we recursively define \emph{the $m$-fold iterated logarithm} by $\log^{(m)}(x) \coloneqq \log\bigl(\log^{(m-1)}(x)\bigr)$ whenever $\log^{(m-1)}(x) \in (1,\infty)$. 
\end{itemize}

Our first theorem quantifies how a near return of the coordinate function to itself, under the semigroup generated by $2$ and $3$, forces the resolution entropy to decay.

\begin{Theorem}[Entropy decay from near return of the coordinate function] \label{thm.inv.return} Let $\mu$ be a $\times 2 \times 3$-invariant Borel probability measure on $\rzee$, let $d \in \enn$ and let $\epsilon \in (0,1]$. Suppose that there exist $j,k \in \{0,\ldots,d\}$, not both zero, such that:
\begin{align} \label{thm.inv.return.conc} \int_{\rzee} \bigl\lvert e(2^j3^kx) - e(x) \bigr\rvert^2 \deee \mu(x) \leq \epsilon \end{align}
Then for all $r \in \ennplus$ satisfying $\log\log(r) \geq (2^{43}/\epsilon)\bigl(3d + 5 + \log_2(1/\epsilon)\bigr)$ we have $\achn(\mu,r) \leq \epsilon$. \end{Theorem}

Theorem \ref{thm.inv.return} will be derived in Section \ref{cha.introapp.four} from Theorem \ref{thm.entret}, which states it in terms of the return defects of Definition \ref{def.ret}.\\

Keeping the notation of Theorem \ref{thm.inv.return}, suppose $\mu$ has discrete spectrum under $x \mapsto 2x$, that is, the $\times 2$-eigenfunctions span a dense subspace of $L^2(\rzee,\mu)$. (See Definition \ref{def.discretespec} and Remark \ref{rem.discretespec} below for details.) Then for every $\epsilon$ there exists $d$ such that (\ref{thm.inv.return.conc}) holds, and so we find $\lim_{r \to \infty}\achn(\mu,r) =0$. This can be interpreted as attaching a return-dependent rate to the Rudolph--Johnson theorem in the discrete spectrum case.\\

The next theorem gives sufficient conditions on proximity of the coordinate function to eigenspaces to obtain entropy rates of decay via Theorem \ref{thm.inv.return}.

\begin{Theorem}[Entropy decay from rate of adherence to eigenfunctions] \label{thm.inv.eigen} Let $\mu$ be a $\times 2 \times 3$-invariant Borel probability measure on $\rzee$.
\begin{enumerate}[label=\upshape(\arabic*)]
\item \label{thm.inv.eigen.c1} Suppose there exists $a \in (0,\infty)$ such that for all but finitely many $n \in \ennplus$ the coordinate function lies within $L^2(\rzee,\mu)$-distance $n^{-a}$ of an $n$-term linear combination of $\times 2$-eigenfunctions. Then we have:
\begin{align*} \limsup_{r \to \infty}\ \frac{\log(\achn(\mu,r))}{\log^{(4)}(r)} \leq -2a \end{align*}
\item \label{thm.inv.eigen.c2} Suppose there exists $c \in (0,\infty)$ such that for all but finitely many $n \in \ennplus$ the coordinate function lies within $L^2(\rzee,\mu)$-distance $e^{-cn}$ of an $n$-term linear combination of $\times 2$-eigenfunctions. Then we have:
\begin{align*} \limsup_{r \to \infty}\ \frac{\log(\achn(\mu,r))}{\sqrt{\log^{(3)}(r)}} \leq -2\sqrt{c} \end{align*}
\end{enumerate} \end{Theorem}

Theorem \ref{thm.inv.eigen} will be derived in Section \ref{cha.introapp.dioph} from Clauses \ref{thm.eigqual.c1} and \ref{thm.eigqual.c2} of Theorem \ref{thm.eigqual}.\\

The roles of $2$ and $3$ here are symmetric: Theorems \ref{thm.inv.eigen}, \ref{thm.inv.torsion} and \ref{thm.inv.atom} hold verbatim with $\times 3$-eigenfunctions in place of $\times 2$-eigenfunctions. Indeed, the same proofs apply to the eigenfunctions of any transformation other than the identity in the $\times 2 \times 3$-semigroup.\\

The next theorem gives an alternative type of sufficient condition for entropy rates of decay via Theorem \ref{thm.inv.return}: we can assume an arbitrary number of eigenfunction terms with controlled eigenvalues. The rate improves on that in Theorem \ref{thm.inv.eigen}.

\begin{Theorem}[Entropy decay from torsion eigenfunctions] \label{thm.inv.torsion} Let $\mu$ be a $\times 2 \times 3$-invariant Borel probability measure on $\rzee$, and suppose there exists $a \in (0,\infty)$ such that for all but finitely many $q \in \ennplus$ there is a positive integer $k \leq q$ for which the coordinate function lies within $L^2(\rzee,\mu)$-distance $q^{-a}$ of a linear combination of $\times 2$-eigenfunctions whose eigenvalues are $k$-th roots of unity. Then we have:
\begin{align*} \limsup_{r \to \infty}\ \frac{\log(\achn(\mu,r))}{\log^{(3)}(r)} \leq -\frac{2a}{2a+1} \end{align*} \end{Theorem}

Theorem \ref{thm.inv.torsion} will be derived in Section \ref{cha.introapp.dioph} from Clause \ref{thm.eigqual.c3} of Theorem \ref{thm.eigqual}. It keeps the same polynomial-adherence hypothesis as Clause \ref{thm.inv.eigen.c1} in Theorem \ref{thm.inv.eigen}, but the torsion hypothesis here removes a logarithm from the rate and nearly reaches the ceiling of the method.\\

Finally, we obtain a bound on the correlation between the coordinate function and any one eigenfunction.

\begin{Theorem}[Entropy caps every eigencomponent of the coordinate function] \label{thm.inv.atom} Let $\mu$ be a $\times 2 \times 3$-invariant Borel probability measure on $\rzee$, let $\epsilon \in (0,1]$ and let $\psi$ be a normalized $\times 2$-eigenfunction. Then for all $r \in \ennplus$ satisfying $\log\log(r) \geq 2^{47}/\epsilon^2$ we have:
\begin{align*} \left\lvert \int_{\rzee} e(x)\,\overline{\psi(x)} \deee \mu(x) \right\rvert^2  \leq 1-\achn(\mu,r) + \epsilon \end{align*} \end{Theorem}

\section{Discussion of AI use}

Here, we repeat the discussion that appears in Section 0.2 of the companion article \cite{RJnodyn}. The artificial intelligence models Claude Opus 4.8, Opus 5, Fable 5 and Fable 5.1 from Anthropic PBC were extensively used in the composition of this article, as follows.
\begin{description}
	\item[\textsf{Mathematical falsification}] Most importantly, to accelerate exclusion of incorrect ideas on the part of the authors.
	\item[\textsf{Mathematical bookkeeping}] Secondarily, to accelerate the accurate management of constants and parameters in proof design.
	\item[\textsf{Latex typography}] Lastly, to accelerate the identification of typographical errors in the Latex manuscript.
\end{description}

\hrulefill
\vspace{12 pt}

The specific formal style in this manuscript is not a marker of AI-produced text. Indeed, quite the opposite is true: the Claude models used cannot reliably reproduce this style. The authors wrote the manuscript line-by-line in order to verify it themselves. Claude was called for additional verification and error-cleaning after every line.\\

The authors are optimistic that this style is a functioning compromise between human readability and AI verifiability. It should be reasonably possible for a graduate student with a course in measure theory to formally check each line to their own satisfaction. It should also be possible for a frontier AI system to verify the manuscript without the need to impute meaning or prove facts independently.

\chapter{Introduction} \label{cha.intro}

\section{Return defects}\label{cha.intro.retdef}

In Section \ref{cha.intro.retdef} we introduce a mechanism to quantify the return of the coordinate function as in (\ref{thm.inv.return.conc}). First, we adopt the following minor notational convenience.

\begin{itemize} \item For $x,y \in \arr \cup\{-\infty,\infty\}$ it will be convenient to adopt the notation $[x,y]_{\zee} \coloneqq \{n \in \zee: x \leq n \leq y \}$. For distinctness we will decorate the standard notations denoting real intervals with a subscript $\mathbb{R}$. Thus, for example, the closed unit interval in the real line will be denoted $[0,1]_{\arr}$ and the strictly positive numbers will be denoted $(0,\infty)_\arr$. \end{itemize}

\begin{Definition}[Return defects] \label{def.ret} For $\mu \in \prob(\rzee)$ and $d \in \ennplus$ we define the \emph{return defect of $\mu$ at box size $d$} to be the element $\ret(\mu \verr d)$ of $[0,2]_\arr$ given as follows.
\begin{align} \label{def.ret.1} \ret(\mu \verr d) \coloneqq \min \Bigl\{ 1 - \re\,\hat{\mu}\bigl(2^{j}3^{k}-1\bigr) : (j,k) \in [0,d]_\zee^2 \setminus \{(0,0)\} \Bigr\} \end{align} \end{Definition}

For $\mu \in \prob(\rzee)$ and $n \in \zee$ we have the following elementary reconfiguration of the quantity at the right of (\ref{def.ret.1}), which illustrates that it is small when $\mu$ concentrates near $\frac{1}{n}\zee$.
\begin{align*} 1 - \re\,\hat{\mu}(n) = 2\int_{\rzee} \sin^2(\pi n x) \deee \mu(x) \end{align*}
Here we identify $\rzee$ with a subset of $\arr$ via the coset transversal $[0,1)_\arr$. We now collect two observations about Definition \ref{def.ret}.

\begin{Remark}[Return defects as distances between exponentials] \label{rem.retdist} Let $\mu \in \prob(\rzee)$ and let $d \in \ennplus$.
\begin{enumerate}[label=\upshape(\arabic*)]
	\item \label{rem.retdist.c1} Elementary rearrangement shows that we have $2(1 - \re\,\hat{\mu}(g-1)) = \mu\bigl(|\mathfrak{e}^{g}-\mathfrak{e}^1|^2\bigr)$ for all $g \in \zee$. Thus $\ret(\mu \verr d) \leq \eta$ asserts that there exists $(j,k) \in [0,d]_\zee^2 \setminus \{(0,0)\}$ satisfying:
	\begin{align*} \mu\bigl(|\mathfrak{e}^{2^{j}3^{k}} - \mathfrak{e}^1|^2\bigr) \leq 2\eta \end{align*}
	\item \label{rem.retdist.c2} Enlarging $d$ enlarges the range of the minimum in (\ref{def.ret.1}), so $\ret(\mu \verr d)$ is nonincreasing in $d$.
	\item \label{rem.retdist.c3} Let $(j,k)$ attain the minimum in (\ref{def.ret.1}) and write $h \coloneqq 2^{j}3^{k}-1$. Then $h \in [1,6^d-1]_\zee$ and $1 - \re\,\hat{\mu}(h) = \ret(\mu \verr d)$, so that $\lvert \hat{\mu}(h) \rvert \geq 1 - \ret(\mu \verr d)$.
\end{enumerate} \end{Remark}

\begin{Remark}[Return defect of Lebesgue measure] \label{rem.lebreturn} For $(j,k) \in \enn^2 \setminus \{(0,0)\}$ we have $2^j3^k - 1 \neq 0$ and so $\hat{\lambda}(2^j3^k-1) = 0$. Thus we find $\ret(\lambda \verr d) = 1$ for all $d \in \ennplus$. \end{Remark}

\section{Resolution entropy}

We now formally introduce the notion of Shannon entropy we will need.

\begin{Definition}[Resolution entropy] \label{def.entropy} Let $n \in \ennplus$.
\begin{itemize} \item For $\ell \in [0,n-1]_\zee$, we define the \emph{index-$\ell$ resolution-$n$ cell in $\rzee$} to be the subset $\mathscr{P}_n(\ell)$ of $\rzee$ corresponding to the subset $[\ell/n,(\ell+1)/n)_{\arr}$ of the coset transversal $[0,1)_\arr$. \end{itemize}
Now, also let $\mu \in \prob(\rzee)$.
\begin{itemize} \item We define the \emph{resolution-$n$ entropy of $\mu$} as follows.
\begin{align}\label{def.entropy.0}\ach(\mu , n) \coloneqq - \sum_{\ell=0}^{n-1} \mu\bigl(\mathscr{P}_n(\ell)\bigr) \log\Bigl(\mu\bigl(\mathscr{P}_n(\ell)\bigr)\Bigr)\end{align}
Here, we understand the conventions that $\log$ denotes the natural logarithm and that $0 \log 0 \coloneqq 0$.
\item We also define the \emph{normalized resolution-$n$ entropy of $\mu$} as follows.
\begin{align} \achn(\mu , n) \coloneqq \begin{dcases}
\hspace{0.5 cm} 1 & \mbox{ if }n=1
\\ \frac{\ach(\mu , n)}{\log(n)} & \mbox{ else }\end{dcases} \label{def.entropy-1} \end{align} \end{itemize} \end{Definition}

We note that it follows from Lemma 2.1.1 and Theorem 2.6.4 of \cite{MR2239987} that $\achn(\mu,r) \in [0,1]_\arr$ for all $\mu \in \prob(\rzee)$ and all $r \in \ennplus$.

\section{Formulation of main theorems} \label{cha.introapp}

\subsection{Entropy decay from $S$-unit near returns} \label{cha.introapp.four}

The following is the basic quantitative relationship between entropy and return defects.

\begin{Theorem}[Entropy decay from a near return at finite resolution] \label{thm.entret} For all $\epsilon \in (0,1]_\arr$, all $d \in \ennplus$, all $\times 2 \times 3$-invariant $\mu \in \prob(\rzee)$ and all $n \in \ennplus$ satisfying:
\begin{align} \label{thm.entret.h1} \log\log(n) \geq \frac{2^{42}}{\epsilon} \log_2\left( \frac{8(6^d+1)}{\epsilon} \right) \end{align}
we have $\achn(\mu,n) \leq \ret(\mu \verr d) + \epsilon$. \end{Theorem}

Theorem \ref{thm.inv.return} follows from Theorem \ref{thm.entret} with $\epsilon$ replaced by $\epsilon/2$, read through Clause \ref{rem.retdist.c1} of Remark \ref{rem.retdist}.\\

\begin{Remark}[One multiplier does not suffice] \label{rem.onemult} The hypothesis of invariance under both $\times 2$ and $\times 3$ in our import Lemma \ref{cor.master} from the companion article \cite{RJnodyn}, and hence in Theorem \ref{thm.entret}, cannot be weakened to invariance under one multiplier. Indeed, for $p \in (0,1/2]_\arr$ let $\nu_p$ be the bias-$p$ ternary Cantor measure, that is, the self-similar measure of Theorem 4.4 in \cite{MR0625600} for the two similitudes of the line with ratio $1/3$ fixing $0$ and $1$, with weights $1-p$ and $p$.\\

This measure is $\times 3$-invariant, and it is the distribution of a point of the middle-thirds Cantor set whose ternary digits are independent and equal to $2$ with probability $p$, by Theorem 4.4(4) of the same reference. Theorem 0.5 of \cite{RJnodyn} gives:
\begin{align*} \lvert \hat{\nu}_p(1) \rvert \geq 1 - \frac{8}{\log(1/p)}\limsup_{r \to \infty}\achn(\nu_p,r)\end{align*}
and the limit superior is positive, since the classical computation of the entropy of a Bernoulli shift (Theorem 2.6.6 of \cite{MR2239987}, transferred to resolution cells as in \cite{RJnodyn}) gives the following for all $k \in \ennplus$.
 \begin{align*} \achn(\nu_p,3^k) = -\frac{p\log(p)+(1-p)\log(1-p)}{\log(3)} \end{align*}
 Applied to $\nu_p$, Lemma \ref{cor.master} would instead give $\lvert \hat{\nu}_p(1) \rvert \leq 1 - \limsup_{r \to \infty}\achn(\nu_p,r)$, which is incompatible with the previous inequality as soon as $p \leq e^{-9}$. \end{Remark}

We now provide the appropriate translations of return defect rates into entropy rates.

\begin{Theorem}[Entropy decay from $S$-unit near returns] \label{thm.retqual} Let $\mu \in \prob(\rzee)$ be $\times 2 \times 3$-invariant.
\begin{enumerate}[label=\upshape(\arabic*)]
\item \label{thm.retqual.c1} We have:
	\begin{align} \label{thm.retqual.h1} \limsup_{r \to \infty}\ \frac{\log\bigl(\achn(\mu,r)\bigr)}{\log^{(4)}(r)}    \leq \limsup_{m \to \infty}\ \frac{\log\bigl(\ret(\mu \verr m)\bigr)}{\log\log(m)} \end{align}
\item \label{thm.retqual.c2} We have:
	\begin{align} \label{thm.retqual.h3} \limsup_{r \to \infty}\ \frac{\log\bigl(\achn(\mu,r)\bigr)}{\sqrt{\log^{(3)}(r)}} \leq \limsup_{m \to \infty}\ \frac{\log\bigl(\ret(\mu \verr m)\bigr)}{\sqrt{\log(m)}} \end{align}
\item \label{thm.retqual.c3} We have:
	\begin{align} \label{thm.retqual.h5} \limsup_{r \to \infty}\ \frac{\log\bigl(\achn(\mu,r)\bigr)}{\log^{(3)}(r)} \leq \Biggl(1 - \limsup_{m \to \infty}\ \frac{\log\bigl(\ret(\mu \verr m)\bigr)}{\log(m)}\Biggr)^{-1} - 1 \end{align}
\end{enumerate} \end{Theorem}

Throughout, we use the conventions $\log(0) \coloneqq -\infty$ and $x/\infty \coloneqq 0$ of \ref{conv.loginf}, so that every display in Theorem \ref{thm.retqual} makes sense when $\ret(\mu \verr m) = 0$ or $\achn(\mu,r) = 0$. In particular the right side of (\ref{thm.retqual.h5}) equals $-1$ when the limit superior there is $-\infty$.

\begin{Remark}[Known invariant measures] \label{rem.lebexclfou} We consider how Theorem \ref{thm.retqual} applies to the known examples of $\times 2 \times 3$-invariant measures.
\begin{itemize} \item If $\mu \in \prob(\rzee)$ is supported on $n$ points then $\achn(\mu,r) \leq \log(n)/\log(r)$ for all $r \in [2,\infty)_\zee$, so the numerators at the left of (\ref{thm.retqual.h1}), (\ref{thm.retqual.h3}) and (\ref{thm.retqual.h5}) are at most $\log\log(n) - \log\log(r)$. Thus the left sides of the inequalities diverge to $-\infty$ and the statements are always satisfied.

\item In contrast, it follows directly from Definition \ref{def.entropy} that we have $\achn(\lambda,r) = 1$ for all $r \in \ennplus$, so that the left side of each inequality in Theorem \ref{thm.retqual} equals $0$ for $\mu = \lambda$. By Remark \ref{rem.lebreturn} we have $\log(\ret(\lambda \verr m)) = 0$ for all $m \in \ennplus$, so that the limit superior at the right of each of (\ref{thm.retqual.h1}), (\ref{thm.retqual.h3}) and (\ref{thm.retqual.h5}) equals $0$. Thus all three inequalities hold with equality for $\mu = \lambda$, which has no entropy decay and has return defect $1$ at every box size. \end{itemize} \end{Remark}

The philosophical reader may contemplate whether understanding Furstenberg's conjecture has proved difficult because one must inherently deal with the exceptionally slow rates that appear in our Theorem \ref{thm.retqual} and other work such as Theorem 1.10 from the article \cite{MR2563089} of Bourgain, Lindenstrauss, Michel and Venkatesh. It is hard to find  much extant analytic theory that involves four-iterate logarithmic rates in a load-bearing way, and what does exist seems heavily concentrated in the theory of primes. (See for instance the bounds on large gaps between consecutive primes due to Rankin \cite{Ran38} and to Ford, Green, Konyagin, Maynard and Tao \cite{FGKMT18}, which carry a four-iterate logarithm, and the oscillation theorem of Littlewood \cite{Lit14} for the distribution of primes, which carries a three-iterate logarithm.)

\subsection{Entropy decay from adherence to eigenfunctions} \label{cha.introapp.dioph}

We now formally introduce the eigenfunctions and associated devices. Again, we emphasize that the present article could be written for eigenfunctions of any nonidentity element of the $\times 2 \times 3$ semigroup.

\begin{Definition}[$\times 2$-eigenfunctions] \label{def.eigen} Let $\mu \in \prob(\rzee)$ be $\times 2$-invariant.
\begin{itemize}
\item For $\theta \in \rzee$ we define the \emph{$\theta$-eigenspace of $\mu$} to be the subspace $\mathcal{E}_\mu(\theta)$ of $L^2(\rzee,\mu)$ consisting of all $\psi \in L^2(\rzee,\mu)$ such that $\psi(2x) = e(\theta)\psi(x)$ for $\mu$-almost every $x \in \rzee$. It is closed, being the kernel of $\eusc{K}_{\mu,2} - e(\theta)\eusc{I}_\mu$ in the notation of Definition \ref{def.koopman}.
\item We define a \emph{$\times 2$-eigenfunction with eigenvalue $\theta$} to be a nonzero element of $\mathcal{E}_\mu(\theta)$, and we define such an eigenfunction to be \emph{normalized} if it is a unit vector of $L^2(\rzee,\mu)$.
\item For $n \in \ennplus$ and a list $\Psi \coloneqq (\psi_1,\ldots,\psi_n)$ of $\times 2$-eigenfunctions, we define the \emph{eigenvalue list of $\Psi$} to be the list $(\theta_1,\ldots,\theta_n)$ in $\rzee$ with $\psi_k \in \mathcal{E}_\mu(\theta_k)$ for all $k \in [1,n]_\zee$, and we define the \emph{eigenresidual of $\mu$ along $\Psi$} as follows, where $\Pi$ denotes orthogonal projection onto the finite-dimensional subspace $\spn(\Psi)$ of $L^2(\rzee,\mu)$.
\begin{align*} \eir(\mu \verr \Psi) \coloneqq \mu\left( \lvertt \mathfrak{e}^1 - \Pi_{\spn(\Psi)} \mathfrak{e}^1 \rvertt^2 \right) \end{align*}
\item For $q \in \ennplus$ we define the \emph{$q$-torsion eigenspace of $\mu$} by $\mathcal{T}_\mu(q) \coloneqq \sum_{p=0}^{q-1} \mathcal{E}_\mu(p/q)$, which is a closed subspace of $L^2(\rzee,\mu)$ by Remark \ref{rem.eigortho}, and we define the \emph{$q$-torsion residual of $\mu$} as follows.
\begin{align*} \tor(\mu \verr q) \coloneqq \min_{1 \leq k \leq q}\ \mu\left( \lvertt \mathfrak{e}^1 - \Pi_{\mathcal{T}_\mu(k)} \mathfrak{e}^1 \rvertt^2 \right) \end{align*}
\end{itemize} \end{Definition}

\begin{Definition}[Discrete spectrum] \label{def.discretespec} Let $\mu \in \prob(\rzee)$ be $\times 2$-invariant. We define $\mu$ to have \emph{discrete spectrum under $\times 2$} if the $\times 2$-eigenfunctions span a dense linear subspace of $L^2(\rzee,\mu)$. \end{Definition}

We will also want a notation for the distance from a point of $\rzee$ to $0$.

\begin{itemize} \item For $x \in \rzee$ we write $\lVert x \rVert_\zee \coloneqq \min_{m \in \zee} \lvert t - m \rvert$, where $t \in \arr$ is any representative of $x$. Thus $\lVert x \rVert_\zee$ is the distance from $x$ to $0$ in the quotient metric of $\rzee$. \end{itemize}

\begin{Remark}[Eigenfunctions of Lebesgue measure] \label{rem.lebexcldioph} Since $\mathtt{T}_2$ is ergodic with respect to $\lambda$, and ergodic endomorphisms of compact groups are mixing, the $\times 2$-eigenfunctions with respect to $\lambda$ are the nonzero constants. (See \S 1.5 and Theorem 1.28 of \cite{MR0648108}.) Thus $\eir(\lambda \verr \Psi) = 1 - \lvert \hat{\lambda}(1) \rvert^2 = 1$ for every list $\Psi$ and $\tor(\lambda \verr q) = 1$ for every $q \in \ennplus$, so that every right side in Theorem \ref{thm.eigqual} below equals $0$ for $\mu = \lambda$. Since also $\achn(\lambda,r) = 1$ for all $r \in \ennplus$, as recorded in Remark \ref{rem.lebexclfou}, the three inequalities of Theorem \ref{thm.eigqual} hold with equality for $\mu = \lambda$, as those of Theorem \ref{thm.retqual} do. \end{Remark}

\begin{Remark}[Discrete spectrum and the Rudolph--Johnson theorem] \label{rem.discretespec} Let $\mu \in \prob(\rzee)$ be $\times 2 \times 3$-invariant with discrete spectrum under $\times 2$. Then for every $\epsilon \in (0,1]_\arr$ there is a finite list $\Psi$ of $\times 2$-eigenfunctions with $\eir(\mu \verr \Psi) \leq \epsilon$, and Lemma \ref{olem.dirichlet} together with Clause \ref{lem.eigret.c1} of Lemma \ref{lem.eigret} then provides $d \in \ennplus$ with $\ret(\mu \verr d) \leq 3\epsilon$. Theorem \ref{thm.entret} therefore gives $\lim_{r \to \infty} \achn(\mu,r) = 0$.\\

Since $\lambda$ does not have discrete spectrum by Remark \ref{rem.lebexcldioph}, this is the conclusion that the Rudolph--Johnson theorem provides for such $\mu$, in the form of vanishing normalized entropy at large resolution. Theorem \ref{thm.eigqual} attaches a rate to it as soon as the approximation of the coordinate function by eigenfunctions is quantified, and this is the sense in which Theorems \ref{thm.inv.eigen} and \ref{thm.inv.torsion} recover the Rudolph--Johnson theorem in the discrete spectrum case. \end{Remark}

Theorem \ref{thm.eigreturn} below is the eigenvalue version of the quantitative return statement in Theorem \ref{thm.entret}.

\begin{Theorem}[Entropy decay from a return of the eigenvalues at finite resolution] \label{thm.eigreturn} For all $\epsilon \in (0,1]_\arr$, all $t \in [0,1/2]_\arr$, all $n,j \in \ennplus$, all $\times 2 \times 3$-invariant $\mu \in \prob(\rzee)$, all lists $\Psi$ of $n$ $\times 2$-eigenfunctions whose eigenvalue list $(\theta_k)_{k=1}^n$ satisfies $\max_{1 \leq k \leq n} \lVert j\theta_k \rVert_\zee \leq t$, and all $r \in \ennplus$ satisfying:
\begin{align} \label{thm.eigreturn.h1} \log\log(r) \geq \frac{2^{42}}{\epsilon} \log_2\left( \frac{8(2^j+1)}{\epsilon} \right) \end{align}
we have $\achn(\mu,r) \leq 2\eir(\mu \verr \Psi) + 2\pi^2t^2 + \epsilon$. \end{Theorem}

We now state the eigenvalue version of the rates in Theorem \ref{thm.retqual}.

\begin{Theorem}[Entropy decay from adherence to eigenfunctions] \label{thm.eigqual} Let $\mu \in \prob(\rzee)$ be $\times 2 \times 3$-invariant, and for $n \in \ennplus$ let $\Psi_n$ be a list of $n$ $\times 2$-eigenfunctions.
\begin{enumerate}[label=\upshape(\arabic*)]
\item \label{thm.eigqual.c1} We have:
\begin{align} \label{thm.eigqual.1} \limsup_{r \to \infty}\ \frac{\log\bigl(\achn(\mu,r)\bigr)}{\log^{(4)}(r)} \leq \limsup_{n \to \infty} \frac{\log(\eir(\mu \verr \Psi_n))}{\log(n)} \end{align}
\item \label{thm.eigqual.c2} We have:
\begin{align} \label{thm.eigqual.2} \limsup_{r \to \infty}\ \frac{\log\bigl(\achn(\mu,r)\bigr)}{\sqrt{\log^{(3)}(r)}} \leq -\sqrt{-2\limsup_{n \to \infty}\ \frac{\log\bigl(\eir(\mu \verr \Psi_n)\bigr)}{n}} \end{align}
\item \label{thm.eigqual.c3} We have:
\begin{align} \label{thm.eigqual.3} \limsup_{r \to \infty}\ \frac{\log\bigl(\achn(\mu,r)\bigr)}{\log^{(3)}(r)} \leq \Biggl(1 - \limsup_{q \to \infty}\ \frac{\log\bigl(\tor(\mu \verr q)\bigr)}{\log(q)}\Biggr)^{-1} - 1 \end{align}
\end{enumerate} \end{Theorem}

It follows directly from Definition \ref{def.eigen} that we have $\eir(\mu \verr \Psi_n) \leq 1$ and $\tor(\mu \verr q) \leq 1$. Hence the limits superior at the right of (\ref{thm.eigqual.2}) and (\ref{thm.eigqual.3}) are at most $0$ and both right sides are well-defined under the conventions of \ref{conv.loginf}.\\

Theorems \ref{thm.inv.eigen} and \ref{thm.inv.torsion} follow from Theorem \ref{thm.eigqual} through the following dictionary. 
\begin{itemize} \item An eigenvalue $\zeta$ in the sense of the Invitation is $e(\theta)$ for the eigenvalue $\theta \in \rzee$ of Definition \ref{def.eigen}, so that the $k$-th roots of unity correspond to the $k$-torsion elements of $\rzee$.
\item The $L^2(\rzee,\mu)$-distance from the coordinate function to a linear combination of eigenfunctions is at least the square root of $\eir(\mu \verr \Psi)$ for the list $\Psi$ of those eigenfunctions, or of $\tor(\mu \verr q)$ when they are $k$-torsion for some $k \in [1,q]_\zee$.
\item Thus a distance bound $n^{-a}$, $e^{-cn}$ or $q^{-a}$ bounds the right side of the corresponding clause of Theorem \ref{thm.eigqual} by $-2a$, $-2\sqrt{c}$ or $-2a/(2a+1)$ respectively, which are the exponents of Theorems \ref{thm.inv.eigen} and \ref{thm.inv.torsion}. \end{itemize}

The proofs use $\times 2$ only through the frequencies $2^j - 1$, which Definition \ref{def.ret} shares with $3^j - 1$ and with $m^j - 1$ for every $m \in \mathbb{S}_{2,3} \setminus \{1\}$; thus Theorems \ref{thm.eigreturn}, \ref{thm.atom} and \ref{thm.eigqual} hold for the eigenfunctions of $\times m$ with the constant $2^R$ in the thresholds of the first two replaced by $m^R$, and Theorem \ref{thm.eigqual} holds verbatim.

\begin{Remark}[Typical and arithmetic eigenvalues] \label{rem.regimes} Clause \ref{thm.eigqual.c2} of Theorem \ref{thm.eigqual} assumes nothing about the eigenvalues of $\Psi_n$, and the only Diophantine input in its proof is Dirichlet's theorem, applied to the whole eigenvalue list of $\Psi_n$. By the theorem of Khinchin and Groshev (see Chapter VII of \cite{MR0087708}), Lebesgue-almost every list of $n$ elements of $\rzee$ admits no simultaneous return to $0$ faster than Dirichlet's theorem provides, up to logarithmic factors.\\

Thus for typical eigenvalue lists the argument for Clause \ref{thm.eigqual.c2} cannot be improved by a better simultaneous return of the whole list. However, this is not a formal claim that the rate in Clause \ref{thm.eigqual.c2} is sharp: a return of the coordinate function needs a good return only of the eigenvalues whose eigencomponents carry most of its mass, and a sharp bound on the return defects need not yield a sharp bound on the entropy.\\

Diophantine conditions of the kind used in KAM theory, which also hold almost everywhere, push in the opposite direction from what the method needs. They bound the simultaneous returns of an eigenvalue list from below at each scale, and so they supply no upper bound on the return defects $\ret(\mu \verr m)$ for Lemma \ref{cor.master} to use. The improvement in Clause \ref{thm.eigqual.c3} comes instead from arithmetic structure, namely the exact returns of torsion eigenvalues.\end{Remark}

We also have the following result comparing entropy to coordinate-eigenfunction correlation. 

\begin{Theorem}[Entropy caps every eigencomponent of the coordinate function] \label{thm.atom} For all $\epsilon \in (0,1]_\arr$, all $R \in \ennplus$, all $\times 2 \times 3$-invariant $\mu \in \prob(\rzee)$, all normalized $\times 2$-eigenfunctions $\psi$ and all $r \in \ennplus$ satisfying:
\begin{align} \label{thm.atom.h1} \log\log(r) \geq \frac{2^{42}}{\epsilon} \log_2\left( \frac{8(2^R+1)}{\epsilon} \right) \end{align}
we have the inequality below.
\begin{align} \label{thm.atom.conc} \lvert \mu(\mathfrak{e}^1\ov{\psi}) \rvert^2 \leq \frac{1}{R+1} + 1 - \achn(\mu,r) + \epsilon \end{align} \end{Theorem}

Theorem \ref{thm.inv.atom} follows from Theorem \ref{thm.atom} with $\epsilon$ replaced by $\epsilon/2$ and $R \coloneqq \lceil 2/\epsilon \rceil$.

\section{Further directions} \label{cha.intro.further}

Since the key threshold in (\ref{cor.master.h1}) below carries a factor $1/\epsilon$, no argument through Lemma \ref{cor.master} can bound $\achn(\mu,r)$ by less than a constant multiple of $\log^{(3)}(r)/\log\log(r)$. We note that the torsion case of Theorem \ref{thm.eigqual} already approaches this rate. The double logarithm enters through the doubly exponential dependence of $N$ on $t$ in Theorem \ref{othm.rigid}.\\

In the sense of Remark \ref{rem.regimes}, the method is limited by simultaneous approximation. For generic eigenvalue lists nothing better than Dirichlet's theorem gives a simultaneous return of the list to $0$, and torsion is the arithmetic structure that gives an exact one. In the intermediate regime one could consider the eigenvalue lists spanning a subgroup of $\rzee$ of finite rank, for which Dirichlet's theorem applied to the generators gives return times polynomial in the precision and rates which are powers of $\log\log(r)$ with exponent depending on the rank. We have not pursued this.\\

The spectral hypotheses are sharp for the method in the sense of Remark \ref{rem.regimes}: for generic eigenvalues nothing better than Dirichlet's theorem returns the coordinate function, and torsion is the arithmetic structure that returns it exactly. Between these lie the eigenvalue lists spanning a subgroup of $\rzee$ of finite rank, for which Dirichlet's theorem applied to the generators gives return times polynomial in the precision and rates which are powers of $\log\log(r)$ with exponent depending on the rank; we have not pursued this.\\

Finally, the arguments use nothing about $2$ and $3$ beyond the irrationality exponent of Lemma \ref{prop.kappa} and the rigidity theorem of \cite{RJnodyn}, which is stated there for every pair of multiplicatively independent multipliers. The theorems of the present article therefore hold for $\times a \times b$-invariant measures with constants depending on $a$ and $b$ through an irrationality exponent for $\log(a)/\log(b)$, which Lemma \ref{olem.bw} supplies in every case.\\

It is also interesting to contemplate generalizations of the current article into more abstract settings. For example, Lindenstrauss's article \cite{MR2827904} places Furstenberg's conjecture in the broader context of homogeneous dynamics, and the article \cite{BP26} of the first author and Jane Panangaden gives a perspective via metabelian groups.

\chapter{Preparatory material} \label{cha.prep}

Chapter \ref{cha.prep} collects what is needed for the proofs in Chapter \ref{cha.proofs} beyond the material already in Chapter \ref{cha.intro}.

\begin{itemize} \item Section \ref{cha.import} fixes conventions and cites the standard facts from linear algebra and measure theory that we rely on.
\item  Section \ref{cha.import.harm.fou} fixes the Fourier-analytic notation and the exponential infimum in which the imported rigidity theorem is phrased.
\item Section \ref{cha.spec} contains the three spectral lemmas that connect eigenfunctions to return defects and to Fourier coefficients, which are the only results of this chapter with proofs.\end{itemize}

\section{Preliminaries} \label{cha.import}

The conventions of Section \ref{cha.import} are standard and are recorded so that the proofs can cite them by label; the reader may prefer to pass directly to Section \ref{cha.import.harm.fou}.

\subsection{Numerical basics}
\label{cha.import.basic}

We fix the following basic conventions throughout the remainder of the article.

\begin{enumerate}[label=\upshape(\thesection:\alph*)]
\item  \label{conv.min} For $x,y \in \arr$ we adopt infix notations for maximum and minimum: 
\begin{align*} x \wedge y &\coloneqq \min(x,y) & x \vee y &\coloneqq \max(x,y) \end{align*}
It will sometimes be convenient to use large prefix versions of these operators. Thus, for example, if $m,n \in \zee$ satisfy $m \leq n$ and $x_m,\ldots,x_n$ is a list of real numbers we define:
\begin{align*} \bigwedge_{j=m}^n x_j & \coloneqq \min  \bigl \{x_j: j \in [m,n]_{\zee}  \bigr \} &  \bigvee_{j=m}^n x_j & \coloneqq \max \bigl \{x_j: j \in [m,n]_{\zee} \bigr \}  \end{align*}
\item  \label{conv.exp} For $x \in \arr$ we adopt the notation $e(x) \coloneqq e^{2 \pi i x}$. For $n \in \zee$, it will be convenient to have a symbol to denote the function $x \mapsto e(nx)$ itself. To this end we define a function $\mathfrak{e}^n: \arr \to \cee$ by $\mathfrak{e}^n(x) \coloneqq e(nx)$. 
\item  \label{conv.period} We will identify $\zee$-periodic functions on $\arr$ (such as $\mathfrak{e}^n$) with the corresponding functions on $\rzee$. We will sometimes identify the group $\rzee$ with its coset transversal $[0,1)_{\arr}$.
\item \label{conv.loginf} We adopt the conventions $\log(0) \coloneqq -\infty$, $x/\infty \coloneqq 0$ for $x \in (0,\infty)_\arr$ and $\sqrt{\infty} \coloneqq \infty$.
\end{enumerate}

\subsection{Imports from linear algebra} \label{cha.import.linalg}

We fix the following conventions regarding linear algebra. Here and in the remainder of the article we adopt the convention that all vector spaces are assumed to have complex number scalars. In the remainder of Section \ref{cha.import.linalg} we fix a vector space $\mathcal{H}$.

\begin{enumerate}[resume,label=\upshape(\thesection:\alph*)]
\item We adopt the notation $0_\mathcal{H}$ for the zero vector in $\mathcal{H}$.
\item  \label{conv.span} For a subset $E$ of $\mathcal{H}$ we define the \emph{span of $E$} to be the linear subspace $\spn_\mathcal{H}(E)$ of $\mathcal{H}$ given as follows.
\begin{align*} \spn_\mathcal{H}(E) \coloneqq \left\{ \sum_{\xi \in F} \alpha(\xi)\xi: F\mbox{ is a finite subset of } E \mbox{ and }\alpha \in \mathbb{C}^F \right\} \end{align*}
The above display is intended to emphasize that we consider only finite linear combinations when constructing spans. If $E = \{\xi_1,\ldots,\xi_n\}$ for some $n \in \ennplus$ and $\xi_1,\ldots,\xi_n \in \mathcal{H}$ we will ease notation by omitting the set brackets and writing $\spn(\xi_1,\ldots,\xi_n)$ for $\spn_\mathcal{H}(E)$.
\item \label{conv.addsub} If $\mathcal{V}$ and $\mathcal{W}$ are subspaces of $\mathcal{H}$, we define the \emph{sum of $\mathcal{V}$ and $\mathcal{W}$} to be the subspace $\mathcal{V}+\mathcal{W}$ of $\mathcal{H}$ given as follows.
\begin{align*}\mathcal{V} + \mathcal{W} \coloneqq \{\xi + \eta: \xi \in \mathcal{V} \mbox{ and }\eta \in \mathcal{W}\}\end{align*}
\item \label{conv.end} We adopt the notation $\operatorname{End}(\mathcal{H})$ to denote the algebra of linear endomorphisms of $\mathcal{H}$. 
\item \label{conv.operator} If $\eusc{A} \in \operatorname{End}(\mathcal{H})$ and $\xi \in \mathcal{H}$ we will typically ease notation by writing $\eusc{A} \xi$ for the element $\eusc{A}(\xi)$ of $\mathcal{H}$. Similarly, if $E \subseteq \mathcal{H}$ we will typically adopt the notation $\eusc{A} E \coloneqq \{\eusc{A}\xi:\xi \in E\}$ to denote the image of $E$ under $\eusc{A}$.
\item \label{conv.zeroid} We adopt the notations $\eusc{I}_\mathcal{H}$ to denote the identity operator on $\mathcal{H}$ and $\eusc{Z}_\mathcal{H}$ to denote the zero operator on $\mathcal{H}$. \end{enumerate}

Let us now assume that $\mathcal{H}$ has the further structure of a Hilbert space. 

\begin{enumerate}[resume,label=\upshape(\thesection:\alph*)]
\item  \label{conv.sesqui} We adopt the sesquilinearity convention that complex inner products are assumed to be linear in the first coordinate and conjugate-linear in the second coordinate.
\item \label{conv.adjoint} If $\mathcal{H}$ is a Hilbert space and $\eusc{A} \in \operatorname{End}(\mathcal{H})$ is bounded, we adopt the notation $\eusc{A}^\ast$ for the adjoint of $\eusc{A}$, defined by the stipulation that $\langle \eusc{A}\xi,\eta \rangle = \langle \xi,\eusc{A}^\ast\eta \rangle$ for all $\xi,\eta \in \mathcal{H}$.
\item \label{conv.perp} If $\mathcal{V}$ is a closed subspace of a Hilbert space $\mathcal{H}$, we adopt the notation $\mathcal{V}^\perp$ for the orthogonal complement of $\mathcal{V}$ and $\Pi_\mathcal{V}$ for the orthogonal projection of $\mathcal{H}$ onto $\mathcal{V}$.
\end{enumerate}

For the Hilbert space material used in the conventions above we refer to \cite{MR1157815}.

\subsection{Imports from elementary analysis} \label{cha.import.harm.ele}

We fix the following conventions regarding metric spaces. To this end, let $X$ be a compact metric space.
\begin{enumerate}[resume,label=\upshape(\thesection:\alph*)]
\item  \label{conv.continuous} We adopt the convention that $C(X)$ denotes the vector space of continuous complex-valued functions on $X$. We endow $C(X)$ with the supremum norm $||\cdot||_{\sup}$ given for $\psi \in C(X)$ by $||\psi||_{\sup} \coloneqq \sup \{|\psi(x)|: x \in X \}$. We also write $C(X)^\dagger$ for the dual of $C(X)$, so an element of $C(X)^\dagger$ is by definition a bounded linear functional on $C(X)$.
\end{enumerate}

We now fix certain conventions regarding measure theory. To this end, let $X$ be a standard Borel space.
\begin{enumerate}[resume,label=\upshape(\thesection:\alph*)]
\item \label{conv.bbr} We adopt the notation $\bbr(X)$ for the vector space of bounded Borel functions from $X$ to $\cee$.
\item  \label{conv.prob} We adopt the notation $\operatorname{FMea}(X)$ to denote the set of Borel measures on $X$ such that $\mu(X) < \infty$, and we adopt the notation $\prob(X)$ to denote the set of all Borel probability measures on $X$. Thus for $\mu \in \operatorname{FMea}(X)$ we have $\mu \in \prob(X)$ if and only if $\mu(X) = 1$.
\item \label{conv.integral} If $\mu \in \operatorname{FMea}(X)$ and $\psi:X \to \cee$ is a $\mu$-integrable function it will sometimes be convenient to adopt the notation $\mu(\psi)$ to denote the complex number corresponding to the integral of $\psi$ over $\mu$.\end{enumerate}
Let us now assume in addition that $X$ also carries a Polish topology.
\begin{enumerate}[resume,label=\upshape(\thesection:\alph*)]
 \item  \label{conv.lp} If $\mu \in \operatorname{FMea}(X)$ and $p \in [1,\infty]_{\arr}$, we adopt the convention that $L^p(X,\mu)$ denotes the complex-valued $L^p$ space of $(X,\mu)$. We emphasize here that the convention for inner products introduced in \ref{conv.sesqui} renders the inner product associated with $\mu$ as $\langle \psi,\phi \rangle \coloneqq \mu(\psi\ov{\phi})$ for $\phi,\psi \in L^2(X,\mu)$. We now simplify some of the notation introduced in Section \ref{cha.import.linalg} for the assignment $\mathcal{H} \coloneqq L^2(X,\mu)$.
 \begin{itemize} \item If $E$ is a subset of $L^2(X,\mu)$ we simplify the notation introduced in \ref{conv.span} above by letting $\spn_\mu(E)$ denote the linear subspace $\spn_{L^2(X,\mu)}(E)$ of $L^2(X,\mu)$ spanned by $E$.
\item We simplify the notation introduced in \ref{conv.zeroid} above by writing $\eusc{I}_\mu$ and $\eusc{Z}_\mu$ respectively for the identity and zero operators on $L^2(X,\mu)$.  \end{itemize}\end{enumerate}

\begin{Remark}[Positivity in $C(X{)}^\dagger$] \label{rem.cdaggerpos} If $X$ is a compact metric space then $\operatorname{FMea}(X)$ can be identified with the set of all positive elements of $C(X)^\dagger$. More explicitly, this means that a functional $\zeta \in C(X)^\dagger$ lies in $\operatorname{FMea}(X)$ if and only if we have $\zeta(\psi) \geq 0$ for all nonnegative continuous functions $\psi: X \to [0,\infty)_{\arr}$. \\

We also note that according to Theorem 7.17 of \cite{MR1681462} the space $C(X)^\dagger$ may be identified with the space of all complex-valued measures on $X$. This means that for all $\zeta \in C(X)^\dagger$ there exist $t_0,\ldots,t_3 \in [0,\infty)_{\arr}$ and $\mu_0,\ldots,\mu_3 \in \prob(X)$ such that $\zeta = \sum_{n=0}^3 i^n\, t_n\,\mu_n$. Chapter 7 of \cite{MR1681462} contains more information about these topics, although for our purposes the object $C(X)^\dagger$ is merely a useful framework to parse expressions such as $\alpha\mu + \beta\nu$ with $\alpha,\beta \in \arr$ and $\mu,\nu \in \prob(X)$. \end{Remark}

\begin{Remark}[Almost-everywhere equivalence] \label{rem.ae} If $X$ is a compact metric space and $\mu \in \operatorname{FMea}(X)$ and $p \in [1,\infty]_{\arr}$ we will typically suppress the distinction between an element $\psi$ of $C(X)$ and its equivalence class in $L^p(X,\mu)$. Nevertheless, this distinction will become relevant in Section \ref{cha.import.harm.fou} below and to handle the issue we define the \emph{$\mu$-almost everywhere quotient map} to be the norm-nonexpanding linear map $\aee_\mu:C(X) \to L^\infty(X,\mu)$ given for $\psi \in C(X)$ as follows.
\begin{align*} \aee_\mu(\psi) \coloneqq \Bigl\{ \phi \in \mathbb{C}^X : \phi \mbox{ is a }\mu\mbox{-measurable function from }X\mbox{ to }\cee \mbox{ with } \mu(|\psi-\phi|) = 0 \Bigr \} \end{align*}   \end{Remark}

We now fix the objects behind the spectral hypotheses of Section \ref{cha.introapp.dioph}: the endomorphisms $\mathtt{T}_n$, the Koopman operator of $\times 2$ on $L^2(\rzee,\mu)$, and the two facts about its eigenspaces that the proofs of Lemmas \ref{lem.eigreturn} and \ref{lem.boxavg} use.

\begin{Definition}[Endomorphism maps] \label{def.invariant} For $n \in \ennplus$ we write $\mathtt{T}_n:\rzee \to \rzee$ for the endomorphism $x \mapsto nx$. \end{Definition}

According to (\ref{prose.intro.1}) we have that $\mu \in \prob(\rzee)$ is $\times n$-invariant if and only if $\mu(\psi \circ \mathtt{T}_n) = \mu(\psi)$ for all $\psi \in C(\rzee)$. 

\begin{Definition}[Koopman operator] \label{def.koopman} Let $\mu \in \prob(\rzee)$ be $\times 2$-invariant. We define the \emph{Koopman operator} $\eusc{K}_{\mu,2} \in \operatorname{End}(L^2(\rzee,\mu))$ by stipulating that $\eusc{K}_{\mu,2}(\psi) \coloneqq \psi \circ \mathtt{T}_2$ for all $\psi \in L^2(\rzee,\mu)$. \end{Definition}

Since $\mu$ is $\times 2$-invariant, we have $\mu(|\psi \circ \mathtt{T}_2|^2) = \mu(|\psi|^2)$ for all $\psi \in L^2(\rzee,\mu)$, so $\eusc{K}_{\mu,2}$ is an isometry of $L^2(\rzee,\mu)$, and polarization gives $\langle \eusc{K}_{\mu,2}\psi,\eusc{K}_{\mu,2}\phi \rangle = \langle \psi,\phi \rangle$ for all $\psi,\phi \in L^2(\rzee,\mu)$.

\begin{Remark}[Eigenspaces of the Koopman operator] \label{rem.eigortho} Let $\mu \in \prob(\rzee)$ be $\times 2$-invariant. We will use the following two standard facts about eigenvectors of the contraction $\eusc{K}_{\mu,2}$ of $L^2(\rzee,\mu)$; see, for example, Corollary 8.7 of \cite{MR3410920}.
\begin{enumerate}[label=\upshape(\arabic*)]
\item \label{rem.eigortho.c1} For all $\theta \in \rzee$ the eigenspace $\mathcal{E}_\mu(\theta)$ is closed, being the kernel of $\eusc{K}_{\mu,2} - e(\theta)\eusc{I}_\mu$, and it satisfies $\eusc{K}_{\mu,2}\mathcal{E}_\mu(\theta) \subseteq \mathcal{E}_\mu(\theta)$. Furthermore, for all $\theta,\varphi \in \rzee$ with $\theta \neq \varphi$ we have $\mathcal{E}_\mu(\theta) \perp \mathcal{E}_\mu(\varphi)$.
\item \label{rem.eigortho.c2} For all $\theta \in \rzee$ and all $\psi \in \mathcal{E}_\mu(\theta)$ we have $\eusc{K}_{\mu,2}^\ast\psi = e(-\theta)\psi$.
\end{enumerate} \end{Remark}

\section{Fourier analysis} \label{cha.import.harm.fou}

Section \ref{cha.import.harm.fou} fixes the Fourier-analytic vocabulary of the article.

\subsection{Fourier analysis: foundations}

\begin{Definition}[Trigonometric polynomials] \label{def.pol} We define the vector space of \emph{trigonometric polynomials} to be the linear subspace $\pol(\rzee)$ of $C(\rzee)$ given as follows.
\begin{align*} \pol(\rzee) \coloneqq \spn_{C(\rzee)}(\{\mathfrak{e}^n:n \in \zee\}) \end{align*} \end{Definition}

Thus $\pol(\rzee)$ is a proper dense linear subspace of $C(\rzee)$. We now introduce some machinery to specify the set of exponentials appearing in a polynomial.
\begin{Definition}[Degree as sets of integers] \label{def.degfinset} Let $\mathscr{F}$ be a subset of $\zee$. We write $\pol(\rzee,\mathscr{F})$ to denote the vector space of \emph{degree-}$\mathscr{F}$ \emph{trigonometric polynomials}, which is defined to be the linear subspace of $C(\rzee)$ given as follows.
\begin{align}\label{eq.poldef}  \pol(\rzee,\mathscr{F}) \coloneqq \spn_{C(\rzee)} \{\mathfrak{e}^n: n \in \mathscr{F} \}   \end{align}
For $n \in \enn$ we will simplify notation by defining $\pol(\rzee,n)  \coloneqq  \pol(\rzee,[0,n]_{\zee})$. \end{Definition}
 
In accordance with Remark \ref{rem.ae}, for $\psi \in \pol(\rzee)$ and $\mu \in \prob(\rzee)$ we will typically regard $\psi$ as an element of $L^2(\rzee,\mu)$ without notational indication, and if such indication is necessary we will use the notation $\aee_\mu(\psi)$.\\

We now stipulate our convention regarding Fourier coefficients over $\rzee$.

\begin{Definition}[Fourier coefficients] \label{def.foucoef} For $n \in \zee$ and $\zeta \in C(\rzee)^\dagger$ we define the \emph{frequency-$n$ Fourier coefficient of $\zeta$} to be the complex number given by $\hat{\zeta}(n) \coloneqq \zeta(\mathfrak{e}^{-n})$. \end{Definition}

\begin{Remark}[Polarization for unimodular exponentials] \label{rem.unipol} Let $\mu \in \prob(\rzee)$.
\begin{enumerate}[label=\upshape(\arabic*)]
	\item \label{rem.unipol.c1} Suppose $\psi,\phi \in \bbr(\rzee)$ satisfy $|\psi(x)| = |\phi(x)| = 1$ for all $x \in \rzee$. We have the elementary identity
	\begin{align*} |\psi - \phi|^2 = |\psi|^2 + |\phi|^2 - 2\re(\psi\ov{\phi}) \end{align*} 
	and integrating the above gives:
\begin{align} \label{rem.unipol.1} \mu\bigl(|\psi-\phi|^2\bigr) = 2 - 2\re\,\mu\bigl(\psi\ov{\phi}\bigr) \end{align}
	\item \label{rem.unipol.c2} Let $j,k \in \zee$ and instantiate Clause \ref{rem.unipol.c1} with the assignment $(\psi,\phi) \coloneqq (\mathfrak{e}^{j},\mathfrak{e}^{k})$. Using the elementary identities $\mathfrak{e}^{j}\ov{\mathfrak{e}^{k}} = \mathfrak{e}^{j-k}$ and $\hat{\mu}(-n) = \ov{\hat{\mu}(n)}$, we find:
\begin{align} \label{rem.unipol.2} \mu\bigl(|\mathfrak{e}^{j}-\mathfrak{e}^{k}|^2\bigr) = 2\bigl(1 - \re\,\hat{\mu}(j-k)\bigr) \end{align}\end{enumerate} \end{Remark}

\label{cha.import.harm.endo}

\begin{Definition}[Endomorphs of measures on $\rzee$] \label{def.endomorph} For $n \in \ennplus$ and $\zeta \in C(\rzee)^\dagger$, we define the \emph{$\times n$-endomorph of $\zeta$} to be the element $n_\circledcirc \zeta$ of $C(\rzee)^\dagger$ given by stipulating that $[n_\circledcirc\zeta](\psi) \coloneqq \zeta(\psi \circ \mathtt{T}_n)$ for all $\psi \in C(\rzee)$. \end{Definition}

Endomorphs are the images of measures under the maps $\mathtt{T}_n$. They appear only in the statement of the imported Theorem \ref{othm.rigid}, and Remark \ref{rem.rigidinv} removes them under invariance.

\begin{Remark}[Endomorphs and change-of-variable formulas] \label{rem.endochangeofvar} Let $n \in \ennplus$ and let $\mu \in \operatorname{FMea}(\rzee)$. 
\begin{enumerate}[label=\upshape(\arabic*)] \item \label{rem.endochangeofvar.c1} By Definition \ref{def.endomorph} the functional $n_\circledcirc \mu$ is positive. Thus by Remark \ref{rem.cdaggerpos} this functional is an element of $\operatorname{FMea}(\rzee)$, and it is the unique such element which satisfies $\widehat{\,n_\circledcirc \mu\,}(\ell) = \hat{\mu}(n\ell)$ for all $\ell \in \zee$.\item \label{rem.endochangeofvar.c2} It is similarly direct that $n_\circledcirc \mu$ is the unique element of $\operatorname{FMea}(\rzee)$ which satisfies the following for all $\psi \in C(\rzee)$.
\begin{align*} \int_{\rzee} \psi(x) \deee [n_\circledcirc\mu](x) = \int_{\rzee} \psi(nx) \deee \mu(x) \end{align*}

Combining the above formula with Definition \ref{def.invariant} we see $\mu$ is $\times n$-invariant if and only if $n_\circledcirc \mu = \mu$.\end{enumerate}\end{Remark}
 
\subsection{Exponential infimum} \label{cha.harm.stat.foudef}

The import from the companion article \cite{RJnodyn} which appears as Theorem \ref{othm.rigid} below will be phrased in terms of a certain infimum of quadratic forms over trigonometric polynomials. 

\begin{Definition}[Exponential infimum] \label{def.foustat} Let $\mu \in \prob(\rzee)$ and let $\mathscr{F}$ be a nonempty finite subset of $\zee$. We define the \emph{bandprofile-$\mathscr{F}$ exponential infimum of $\mu$} as follows.
\begin{align*} \einf(\mu , \mathscr{F}) \coloneqq \inf \left \{ \mu \left ( \left \vert \sum_{j \in \mathscr{F}} \alpha(j) \,\mathfrak{e}^j \right \vert^2  \right): \alpha \in \cee^\mathscr{F} \mbox{ and } \sum_{j \in \mathscr{F}} |\alpha(j)|^2 = 1  \right \} \end{align*}
We simplify notation when $\mathscr{F} = [0,n]_{\zee}$ for some $n \in \enn$ by writing the input to the second coordinate simply as $n$. Thus, for example, we will write $\einf(\mu,n)$ instead of $\einf(\mu,[0,n]_{\zee})$. \end{Definition}

We now collect two observations about the above quantity.

\begin{Remark}[Alternative expressions for the exponential infimum] \label{rem.einf} Let $\mu \in \prob(\rzee)$ and let $\mathscr{F}$ be a nonempty finite subset of $\zee$. 
\begin{enumerate}[label=\upshape(\arabic*)]
	\item \label{rem.einf.c1} The exponential infimum may be expressed more directly in terms of Fourier coefficients as follows.
\begin{align*} \einf(\mu , \mathscr{F}) =\inf \left \{ \sum_{j,k \in \mathscr{F}} \alpha(j) \ov{\alpha(k)}\,\hat{\mu}(k-j) :\alpha \in \cee^\mathscr{F} \mbox{ and } \sum_{j \in \mathscr{F}} |\alpha(j)|^2 = 1 \right \} \end{align*}
\item \label{rem.einf.c2} Using the elementary Plancherel identity we have a further alternative expression for the exponential infimum:
\begin{align*} \einf(\mu , \mathscr{F}) = \inf \left \{ \mu \left ( \left \vert \sum_{j \in \mathscr{F}} \alpha(j) \,\mathfrak{e}^j \right \vert^2  \right): \alpha \in \cee^\mathscr{F} \mbox{ and }  \lambda \left ( \left \vert \sum_{j \in \mathscr{F}} \alpha(j) \,\mathfrak{e}^j \right \vert^2  \right) =1 \right \}  \end{align*}
We incorporate the notation of Definition \ref{def.degfinset} by rewriting the last display as follows.
\begin{align} \einf(\mu , \mathscr{F}) = \inf \left \{ \mu(\lvertt \psi \rvertt^2): \psi \in \pol(\rzee,\mathscr{F}) \mbox{ and } \lambda(|\psi|^2) =1 \right \} \label{rem.einf.1} \end{align}
\end{enumerate}\end{Remark}

Remark \ref{rem.foueinf} converts a lower bound on $\einf$ into an upper bound on a single Fourier coefficient, which is the form used throughout the article.

\begin{Remark}[Fourier coefficients bounded by the exponential infimum] \label{rem.foueinf} Let $\mu \in \prob(\rzee)$, let $r \in \ennplus$ and let $h \in [1,r]_\zee$. Choose $\theta \in \arr$ with $e(-\theta)\,\hat{\mu}(h) = |\hat{\mu}(h)|$, and define $\alpha \in \cee^{[0,r]_\zee}$ by $\alpha(0) \coloneqq 1/\sqrt{2}$, $\alpha(h) \coloneqq -e(\theta)/\sqrt{2}$ and $\alpha(j) \coloneqq 0$ for $j \in [1,r]_\zee \setminus \{h\}$. Then $\sum_{j=0}^{r}|\alpha(j)|^2 = 1$, and the sum in Clause \ref{rem.einf.c1} of Remark \ref{rem.einf} evaluates as follows.
\begin{align*} \sum_{j,k=0}^{r} \alpha(j)\ov{\alpha(k)}\,\hat{\mu}(k-j) = \frac{1}{2} + \frac{1}{2} - \frac{1}{2}e(-\theta)\,\hat{\mu}(h) - \frac{1}{2}e(\theta)\,\ov{\hat{\mu}(h)} = 1 - |\hat{\mu}(h)| \end{align*}
Thus Clause \ref{rem.einf.c1} of Remark \ref{rem.einf} gives $\einf(\mu,[0,r]_\zee) \leq 1 - |\hat{\mu}(h)|$, and since $h$ was arbitrary we obtain:
\begin{align} \label{rem.foueinf.1} \max_{1 \leq h \leq r} |\hat{\mu}(h)| \leq 1 - \einf(\mu,[0,r]_\zee) \end{align} \end{Remark}

\section{Spectral theory} \label{cha.spec}

Section \ref{cha.spec} contains the three lemmas that turn spectral information into Fourier information.
\begin{itemize} \item  Lemma \ref{lem.eigreturn} shows that a simultaneous near return of finitely many eigenvalues forces a large real part of a Fourier coefficient at a frequency of the form $2^j - 1$, up to the distance of the coordinate function from the corresponding eigenspaces.
\item  Lemma \ref{lem.eigret} restates this as a bound on return defects.
\item Lemma \ref{lem.boxavg} bounds the correlation of the coordinate function with a single eigenfunction by an average of Fourier coefficients. \end{itemize} 
The first two lemmas above feed the proofs of the theorems of Section \ref{cha.introapp.dioph}, and the third feeds Theorem \ref{thm.atom}.

\subsection{Simultaneous Diophantine approximation}

The following is Dirichlet's theorem on simultaneous approximation, in the form delivered by the box argument of \S 11.12 of \cite{MR2445243}. It is the only Diophantine input of the article, and it supplies the return times in the proof of Theorem \ref{thm.eigqual}.

\begin{oLemma}[Dirichlet's theorem on simultaneous approximation] \label{olem.dirichlet} Let $n,q \in \ennplus$ and let $(\theta_k)_{k=1}^n$ be a list in $\rzee$. Then there exists $j \in [1,q^n]_\zee$ such that:
\begin{align} \label{olem.dirichlet.conc} \bigvee_{k=1}^{n} \lVert j\theta_k \rVert_\zee \leq \frac{1}{q} \end{align} \end{oLemma}

\subsection{Return of the coordinate function from an almost invariant subspace}

Lemma \ref{lem.eigreturn} is the mechanism behind every spectral hypothesis in the article. Keeping the notation of the proof below, we find that if the eigenvalues in a finite set $\Theta$ all nearly return to $0$ under multiplication by $j$ then $\eusc{K}_{\mu,2}^j$ nearly fixes the sum of the corresponding eigenspaces. Thus the coordinate function, if it lies close to that sum, nearly returns to itself along $x \mapsto 2^jx$. Theorem \ref{thm.eigreturn} and Lemma \ref{lem.eigret} are the two consumers of Lemma \ref{lem.eigreturn}.

\begin{Lemma}[Return of the coordinate function from an almost invariant subspace] \label{lem.eigreturn} Let $\mu \in \prob(\rzee)$ be $\times 2$-invariant, let $j \in \ennplus$, let $t \in [0,1/2]_\arr$, let $\Theta$ be a finite subset of $\rzee$ with $\lVert j\theta \rVert_\zee \leq t$ for all $\theta \in \Theta$, and let $\mathcal{V}$ be a closed subspace of $\sum_{\theta \in \Theta} \mathcal{E}_\mu(\theta)$.
Then we have:
\begin{align} \label{lem.eigreturn.conc} \re\,\hat{\mu}(2^j-1) \geq 1 - 2\,\mu\left( \lvertt \mathfrak{e}^1 - \Pi_\mathcal{V}(\mathfrak{e}^1) \rvertt^2 \right) - 2\pi^2t^2 \end{align} \end{Lemma}

The proof below uses Definition \ref{def.koopman} and Clauses \ref{rem.eigortho.c1} and \ref{rem.eigortho.c2} of Remark \ref{rem.eigortho}.

\begin{Proof}[of Lemma \ref{lem.eigreturn}] Let $\mu,j,t,\Theta$ and $\mathcal{V}$ be as in the statement of Lemma \ref{lem.eigreturn}, write $\eusc{K} \coloneqq \eusc{K}_{\mu,2}$, and regard $\mathfrak{e}^1$ as an element of $L^2(\rzee,\mu)$. Also write $\mathcal{W} \coloneqq \sum_{\theta \in \Theta} \mathcal{E}_\mu(\theta)$. By Clause \ref{rem.eigortho.c1} of Remark \ref{rem.eigortho}, $\mathcal{W}$ is a finite orthogonal sum of closed subspaces, hence closed, and $\eusc{K}\mathcal{W} \subseteq \mathcal{W}$. \\

Adopt the notations $\xi \coloneqq \Pi_\mathcal{W}\mathfrak{e}^1$ and $\eta \coloneqq \mathfrak{e}^1 - \xi$. Also write $\xi = \sum_{\theta \in \Theta} \xi_\theta$ with $\xi_\theta \in \mathcal{E}_\mu(\theta)$. Using the Plancherel identity we find:
\begin{align} \label{lem.eigreturn.p1} ||\eta||^2 = 1 - \sum_{\theta \in \Theta} ||\xi_\theta||^2\end{align}
Since $\mathcal{V} \subseteq \mathcal{W}$ we have:
\begin{align} \label{lem.eigreturn.p2} ||\eta||^2 \leq \mu\left( \lvertt \mathfrak{e}^1 - \Pi_\mathcal{V}(\mathfrak{e}^1) \rvertt^2 \right) \end{align}
It follows directly from the definition of the Koopman operator in Definition \ref{def.koopman} that $\eusc{K}^j\mathfrak{e}^1 = \mathfrak{e}^{2^j}$, and so we find:
\begin{align} \label{lem.eigreturn.p3} \re\,\hat{\mu}(2^j-1) =  \re\,\langle \eusc{K}^j\mathfrak{e}^1,\mathfrak{e}^1 \rangle \end{align}
We next record the invariance of $\mathcal{W}^\perp$. Let $\zeta \in \mathcal{W}^\perp$ and $w \in \mathcal{W}$, and write $w = \sum_{\theta \in \Theta} w_\theta$ with $w_\theta \in \mathcal{E}_\mu(\theta)$. By Clause \ref{rem.eigortho.c2} of Remark \ref{rem.eigortho} we have:
\begin{align} \label{lem.eigreturn.p4} \langle \eusc{K}\zeta,w \rangle = \sum_{\theta \in \Theta} \langle \zeta,\eusc{K}^\ast w_\theta \rangle = \sum_{\theta \in \Theta} e(\theta)\langle \zeta,w_\theta \rangle = 0 \end{align}
Thus $\eusc{K}\mathcal{W}^\perp \subseteq \mathcal{W}^\perp$. Iterating the two invariances, $\eusc{K}^j\xi \in \mathcal{W}$ and $\eusc{K}^j\eta \in \mathcal{W}^\perp$, so that $\langle \eusc{K}^j\xi,\eta \rangle = 0 = \langle \eusc{K}^j\eta,\xi \rangle$, and expanding $\mathfrak{e}^1 = \xi + \eta$ in (\ref{lem.eigreturn.p3}) we obtain:
\begin{align} \label{lem.eigreturn.p5} \re\,\hat{\mu}(2^j-1) = \re\,\langle \eusc{K}^j\xi,\xi \rangle + \re\,\langle \eusc{K}^j\eta,\eta \rangle \end{align}
We bound the two terms separately. Since $\eusc{K}$ is an isometry, the Cauchy-Schwarz inequality gives:
\begin{align} \label{lem.eigreturn.p6} \re\,\langle \eusc{K}^j\eta,\eta \rangle \geq -||\eusc{K}^j\eta||\,||\eta|| = -||\eta||^2 \end{align}
By Clause \ref{rem.eigortho.c1} of Remark \ref{rem.eigortho} the summands are pairwise orthogonal, and $\eusc{K}^j\xi_\theta = e(j\theta)\xi_\theta$. We compute:
\begin{align} \re\,\langle \eusc{K}^j\xi,\xi \rangle & = \sum_{\theta \in \Theta} \re\bigl(e(j\theta)\bigr)\,||\xi_\theta||^2 \label{lem.eigreturn.p7.1}
\\ & \geq \bigl(1 - 2\pi^2t^2\bigr) \sum_{\theta \in \Theta} ||\xi_\theta||^2 \label{lem.eigreturn.p7.2}
\\ & = \bigl(1 - 2\pi^2t^2\bigr)\bigl(1 - ||\eta||^2\bigr) \label{lem.eigreturn.p7.3}
\\ & \geq 1 - ||\eta||^2 - 2\pi^2t^2 \label{lem.eigreturn.p7} \end{align}
This computation may be justified as follows.
\begin{itemize}
\item The equality in (\ref{lem.eigreturn.p7.1}) holds by expanding the inner product, and noting that the cross terms vanish by orthogonality.
\item (\ref{lem.eigreturn.p7.2}) follows from (\ref{lem.eigreturn.p7.1}) since for $\theta \in \Theta$ we can use the inequality $\lvert \sin(\pi x)\rvert \leq \pi\lVert x \rVert_\zee$ to find:
\begin{align*} \re\,e(j\theta) = 1 - 2\sin^2(\pi j\theta) \geq 1 - 2\pi^2t^2 \end{align*}
\item (\ref{lem.eigreturn.p7.3}) follows from (\ref{lem.eigreturn.p7.2}) by (\ref{lem.eigreturn.p1}).
\end{itemize}
We compute again:
\begin{align} \re\,\hat{\mu}(2^j-1) & = \re\,\langle \eusc{K}^j\xi,\xi \rangle + \re\,\langle \eusc{K}^j\eta,\eta \rangle \label{lem.eigreturn.p8.1}
\\ & \geq 1 - 2||\eta||^2 - 2\pi^2t^2 \label{lem.eigreturn.p8.3}
\\ & \geq 1 - 2\,\mu\left(\lvertt \mathfrak{e}^1 - \Pi_\mathcal{V}(\mathfrak{e}^1) \rvertt^2\right) - 2\pi^2t^2 \label{lem.eigreturn.p8} \end{align}
This computation may be justified as follows.
\begin{itemize}
\item The equality in (\ref{lem.eigreturn.p8.1}) holds by (\ref{lem.eigreturn.p5}).
\item (\ref{lem.eigreturn.p8.3}) follows from (\ref{lem.eigreturn.p8.1}) by (\ref{lem.eigreturn.p7}) and (\ref{lem.eigreturn.p6}).
\item (\ref{lem.eigreturn.p8}) follows from (\ref{lem.eigreturn.p8.3}) by (\ref{lem.eigreturn.p2}).
\end{itemize}
The computation ending in (\ref{lem.eigreturn.p8}) completes the proof of Lemma \ref{lem.eigreturn}. \end{Proof}

\subsection{Return defects from adherence to eigenfunctions}

Lemma \ref{lem.eigret} turns the Fourier bound of Lemma \ref{lem.eigreturn} into a consequence for the return defects of Definition \ref{def.ret}, which is the form consumed by the three certificates in the proof of Theorem \ref{thm.eigqual}. Clause \ref{lem.eigret.c2} is the torsion case, where the return is exact.

\begin{Lemma}[Return defects from adherence to eigenfunctions] \label{lem.eigret} Let $\mu \in \prob(\rzee)$ be $\times 2$-invariant.
\begin{enumerate}[label=\upshape(\arabic*)]
\item \label{lem.eigret.c1} Let $n,j \in \ennplus$, let $t \in [0,1/2]_\arr$ and let $\Psi$ be a list of $n$ $\times 2$-eigenfunctions whose eigenvalue list $(\theta_k)_{k=1}^n$ satisfies $\max_{1 \leq k \leq n} \lVert j\theta_k \rVert_\zee \leq t$. Then for all $m \in [j,\infty)_\zee$ we have:
\begin{align} \label{lem.eigret.conc1} \ret(\mu \verr m) \leq 2\eir(\mu \verr \Psi) + 2\pi^2t^2 \end{align}
\item \label{lem.eigret.c2} Let $q \in \ennplus$. Then for all $m \in [q,\infty)_\zee$ we have:
\begin{align} \label{lem.eigret.conc2} \ret(\mu \verr m) \leq 2\tor(\mu \verr q) \end{align}
\end{enumerate} \end{Lemma}

The proof below uses Definition \ref{def.ret}, Definition \ref{def.eigen} and Lemma \ref{lem.eigreturn}.

\begin{Proof}[of Lemma \ref{lem.eigret}] Let $\mu$ be as in the statement of Lemma \ref{lem.eigret}.\\

\underline{Verification of Clause \ref{lem.eigret.c1}}:\\

Let $n,j,t$ and $\Psi$ be as in Clause \ref{lem.eigret.c1}, and write $\Theta \coloneqq \{\theta_1,\ldots,\theta_n\}$ and $\mathcal{V} \coloneqq \spn(\Psi)$. We note that $\mathcal{V}$ is a finite-dimensional subspace of $\sum_{\theta \in \Theta} \mathcal{E}_\mu(\theta)$, and that we have $\lVert j\theta \rVert_\zee \leq t$ for all $\theta \in \Theta$ by hypothesis.\\

Thus Lemma \ref{lem.eigreturn} applies, and since $\mu(|\mathfrak{e}^1 - \Pi_\mathcal{V}\mathfrak{e}^1|^2) = \eir(\mu \verr \Psi)$ by Definition \ref{def.eigen}, the conclusion (\ref{lem.eigreturn.conc}) reads:
\begin{align} \label{lem.eigret.p1} \re\,\hat{\mu}(2^j-1) \geq 1 - 2\eir(\mu \verr \Psi) - 2\pi^2t^2 \end{align}
Now let $m \in [j,\infty)_\zee$ and compute:
\begin{align} \ret(\mu \verr m) & \leq 1 - \re\,\hat{\mu}(2^j-1) \label{lem.eigret.p2.1}
\\ & \leq 2\eir(\mu \verr \Psi) + 2\pi^2t^2 \label{lem.eigret.p2} \end{align}
This computation may be justified as follows.
\begin{itemize}
\item The inequality in (\ref{lem.eigret.p2.1}) holds by (\ref{def.ret.1}), since $1 \leq j \leq m$ gives $(j,0) \in [0,m]_\zee^2 \setminus \{(0,0)\}$ and $2^j3^0 - 1 = 2^j - 1$.
\item (\ref{lem.eigret.p2}) follows from (\ref{lem.eigret.p2.1}) by (\ref{lem.eigret.p1}).
\end{itemize}
The inequality between the outer terms of the computation ending in (\ref{lem.eigret.p2}) is (\ref{lem.eigret.conc1}), which completes the verification of Clause \ref{lem.eigret.c1}.\\

\underline{Verification of Clause \ref{lem.eigret.c2}}:\\

Let $q$ be as in Clause \ref{lem.eigret.c2}. By Definition \ref{def.eigen} we may choose $k \in [1,q]_\zee$ such that $\mu(|\mathfrak{e}^1 - \Pi_{\mathcal{T}_\mu(k)}\mathfrak{e}^1|^2) = \tor(\mu \verr q)$, and we write $\Theta \coloneqq \{p/k : p \in [0,k-1]_\zee\}$ and $\mathcal{V} \coloneqq \mathcal{T}_\mu(k)$. By Definition \ref{def.eigen} we have $\mathcal{V} = \sum_{\theta \in \Theta} \mathcal{E}_\mu(\theta)$, and so $\mathcal{V}$ is closed. Moreover $\lVert k\theta \rVert_\zee = 0$ for all $\theta \in \Theta$, since $k \cdot (p/k) = p \in \zee$.\\

Thus Lemma \ref{lem.eigreturn} applies with the assignment $(j,t) \coloneqq (k,0)$, and since $\mu(|\mathfrak{e}^1 - \Pi_{\mathcal{T}_\mu(k)}\mathfrak{e}^1|^2) = \tor(\mu \verr q)$ by the choice of $k$, the conclusion (\ref{lem.eigreturn.conc}) reads:
\begin{align} \label{lem.eigret.p3} \re\,\hat{\mu}(2^k-1) \geq 1 - 2\tor(\mu \verr q) \end{align}
Now let $m \in [q,\infty)_\zee$ and compute:
\begin{align} \ret(\mu \verr m) & \leq 1 - \re\,\hat{\mu}(2^k-1) \label{lem.eigret.p4.1}
\\ & \leq 2\tor(\mu \verr q) \label{lem.eigret.p4} \end{align}
This computation may be justified as follows.
\begin{itemize}
\item The inequality in (\ref{lem.eigret.p4.1}) holds by (\ref{def.ret.1}), since $1 \leq k \leq q \leq m$ gives $(k,0) \in [0,m]_\zee^2 \setminus \{(0,0)\}$.
\item (\ref{lem.eigret.p4}) follows from (\ref{lem.eigret.p4.1}) by (\ref{lem.eigret.p3}).
\end{itemize}
The inequality between the outer terms of the computation ending in (\ref{lem.eigret.p4}) is (\ref{lem.eigret.conc2}). This completes the verification of the two clauses and thereby the proof of Lemma \ref{lem.eigret}. \end{Proof}

\subsection{Eigencomponents of the coordinate function are bounded by a box average}

Lemma \ref{lem.boxavg} is independent of the return mechanism. It bounds the squared correlation of the coordinate function with a normalized eigenfunction by the reciprocal of a box length plus the largest Fourier coefficient in that box, and it is what Theorem \ref{thm.atom} adds to Lemma \ref{cor.master}.

\begin{Lemma}[Eigencomponents of the coordinate function are bounded by a box average] \label{lem.boxavg} Let $\mu \in \prob(\rzee)$ be $\times 2$-invariant, let $R \in \ennplus$ and let $\psi$ be a normalized $\times 2$-eigenfunction. Then we have:
\begin{align} \label{lem.boxavg.conc} \lvert \langle \mathfrak{e}^1,\psi \rangle \rvert^2 \leq \frac{1}{R+1} + \bigvee_{h=1}^{2^R-1} \lvert \hat{\mu}(h) \rvert \end{align} \end{Lemma}

The proof below uses Clause \ref{rem.eigortho.c2} of Remark \ref{rem.eigortho} and Definition \ref{def.foucoef}.

\begin{Proof}[of Lemma \ref{lem.boxavg}] Let $\mu,R$ and $\psi$ be as in the statement of Lemma \ref{lem.boxavg}, let $\theta \in \rzee$ be the eigenvalue of $\psi$, and write $\eusc{K} \coloneqq \eusc{K}_{\mu,2}$. Also let $k \in [0,R]_\zee$.\\

We compute:
\begin{align} \langle \mathfrak{e}^{2^k},\psi \rangle & = \langle \eusc{K}^k\mathfrak{e}^1,\psi \rangle \label{lem.boxavg.p1.1}
\\ & = \langle \mathfrak{e}^1,(\eusc{K}^\ast)^k\psi \rangle \label{lem.boxavg.p1.2}
\\ & = e(k\theta)\,\langle \mathfrak{e}^1,\psi \rangle \label{lem.boxavg.p1} \end{align}
This computation may be justified as follows.
\begin{itemize}
\item The equality in (\ref{lem.boxavg.p1.1}) holds since $\eusc{K}^k\mathfrak{e}^1 = \mathfrak{e}^{2^k}$.
\item (\ref{lem.boxavg.p1.2}) follows from (\ref{lem.boxavg.p1.1}) by the adjoint identity in \ref{conv.adjoint}.
\item (\ref{lem.boxavg.p1}) follows from (\ref{lem.boxavg.p1.2}) since Clause \ref{rem.eigortho.c2} of Remark \ref{rem.eigortho} gives  $(\eusc{K}^\ast)^k\psi = e(-k\theta)\psi$.
\end{itemize}
We define $\Phi \in L^2(\rzee,\mu)$ as follows.
\begin{align} \label{lem.boxavg.p2} \Phi \coloneqq \frac{1}{R+1}\sum_{k=0}^{R} e(-k\theta)\,\mathfrak{e}^{2^k} \end{align}
Multiplying (\ref{lem.boxavg.p1}) by $e(-k\theta)$ and averaging over $k \in [0,R]_\zee$, we obtain:
\begin{align} \label{lem.boxavg.p3} \langle \Phi,\psi \rangle = \langle \mathfrak{e}^1,\psi \rangle \end{align}
We compute:
\begin{align} \lvert \langle \mathfrak{e}^1,\psi \rangle \rvert^2 & = \lvert \langle \Phi,\psi \rangle \rvert^2 \label{lem.boxavg.p4.1}
\\ & \leq ||\Phi||^2 \label{lem.boxavg.p4.2}
\\ & = \frac{1}{(R+1)^2}\sum_{k,\ell=0}^{R} e\bigl((\ell-k)\theta\bigr)\,\hat{\mu}\bigl(2^\ell - 2^k\bigr) \label{lem.boxavg.p4.3}
\\ & \leq \frac{1}{R+1} + \frac{1}{(R+1)^2}\sum_{k,\ell \in [0,R]_\zee : k \neq \ell} \lvert \hat{\mu}(2^\ell - 2^k) \rvert \label{lem.boxavg.p4.4}
\\ & \leq \frac{1}{R+1} + \bigvee_{h=1}^{2^R-1} \lvert \hat{\mu}(h) \rvert \label{lem.boxavg.p4} \end{align}
This computation may be justified as follows.
\begin{itemize}
\item The equality in (\ref{lem.boxavg.p4.1}) holds by (\ref{lem.boxavg.p3}).
\item (\ref{lem.boxavg.p4.2}) follows from (\ref{lem.boxavg.p4.1}) by the Cauchy-Schwarz inequality, since $||\psi|| = 1$.
\item (\ref{lem.boxavg.p4.3}) follows from (\ref{lem.boxavg.p4.2}) by expanding $||\Phi||^2$ from (\ref{lem.boxavg.p2}), using $\langle \mathfrak{e}^{2^k},\mathfrak{e}^{2^\ell} \rangle = \mu(\mathfrak{e}^{2^k-2^\ell}) = \hat{\mu}(2^\ell-2^k)$ by Definition \ref{def.foucoef}.
\item (\ref{lem.boxavg.p4.4}) follows from (\ref{lem.boxavg.p4.3}) since the $R+1$ diagonal terms equal $\hat{\mu}(0) = 1$ and the off-diagonal terms have modulus at most $\lvert \hat{\mu}(2^\ell-2^k) \rvert$.
\item (\ref{lem.boxavg.p4}) follows from (\ref{lem.boxavg.p4.4}) since $\lvert 2^\ell - 2^k \rvert \in [1,2^R-1]_\zee$ for $k \neq \ell$ in $[0,R]_\zee$.\end{itemize}
The inequality between the outer terms of the computation ending in (\ref{lem.boxavg.p4}) is (\ref{lem.boxavg.conc}), which completes the proof of Lemma \ref{lem.boxavg}. \end{Proof}

\chapter{Proofs of theorems} \label{cha.proofs}

\section{Infrastructural results} \label{cha.rigid}

Section \ref{cha.rigid} assembles the input from \cite{RJnodyn} in the form the proofs use. 
\begin{itemize} \item Lemma \ref{prop.kappa} supplies the irrationality exponent on which the constants of the imported Theorem \ref{othm.rigid} depend.
\item  Lemma \ref{cor.master} converts that theorem into a bound on Fourier coefficients.
\item Lemmas \ref{lem.scales} and \ref{lem.engine} are the two devices through which the qualitative theorems are proved.\end{itemize}

\subsection{Obtaining an irrationality exponent}

The constants of Theorem \ref{othm.rigid} below are expressed through an irrationality exponent for $\log(2)/\log(3)$, and Lemma \ref{prop.kappa} produces one. The value $2^{31}$ is convenient rather than optimal. \\

We obtain such an exponent from the theorem of A.~Baker and G.~W\"ustholz on p.~20 of \cite{MR1234835}, which we import in the case of two logarithms of rational integers, where the field has degree one and the heights appearing there reduce to the expressions below.

\begin{oLemma}[Baker-W\"ustholz for two logarithms of integers] \label{olem.bw} Let $a_1,a_2 \in [2,\infty)_\zee$ and let $b_1,b_2 \in \zee$ satisfy $b_1\log(a_1) + b_2\log(a_2) \neq 0$. Then we have the inequality below, where $C \coloneqq 18 \cdot 3! \cdot 2^3 \cdot 32^4 \log(4)$.
\begin{align} \label{olem.bw.conc} \begin{split} & \log\bigl\lvert b_1\log(a_1) + b_2\log(a_2) \bigr\rvert 
\\ & \hspace{1 in} \geq -C \max\{\log(a_1),1\} \max\{\log(a_2),1\} \max\bigl\{\log\bigl(\max\{|b_1|,|b_2|\}\bigr),1\bigr\} \end{split} \end{align} \end{oLemma}

We now state our specialization of the above lemma.

\begin{Lemma}[An irrationality exponent for {$\log(2)/\log(3)$}] \label{prop.kappa} For all $p \in \zee$ and all $q \in [2,\infty)_\zee$ we have:
\begin{align} \label{prop.kappa.conc} \left\lvert \frac{\log(2)}{\log(3)} - \frac{p}{q} \right\rvert \geq q^{-2^{31}} \end{align} \end{Lemma}

Thus $\kappa \coloneqq 2^{31}$ is admissible above, and $8\kappa = 2^{34}$ is the constant in (\ref{othm.rigid.h1}). The proof below uses Lemma \ref{olem.bw}.

\begin{Proof}[of Lemma \ref{prop.kappa}] Let $p$ and $q$ be as in the statement of Lemma \ref{prop.kappa}.\\

\underline{The case $|p| \leq q$}:\\

Since $2$ and $3$ are multiplicatively independent, $q\log(2) - p\log(3) \neq 0$, so Lemma \ref{olem.bw} applies with the assignment $(a_1,a_2,b_1,b_2) \coloneqq (2,3,q,-p)$. We compute:
\begin{align} \log\bigl\lvert q\log(2) - p\log(3) \bigr\rvert & \geq -C\log(3)\max\{\log(q),1\} \label{prop.kappa.p1.1}
\\ & \geq -\frac{C\log(3)}{\log(2)}\log(q) \label{prop.kappa.p1.2}
\\ & = -27 \cdot 2^{26}\log(3)\log(q) \label{prop.kappa.p1} \end{align}
This computation may be justified as follows.
\begin{itemize}
\item The inequality in (\ref{prop.kappa.p1.1}) holds by (\ref{olem.bw.conc}), since $\max\{\log(2),1\} = 1$, $\max\{\log(3),1\} = \log(3)$ and $\max\{|q|,|p|\} = q$.
\item (\ref{prop.kappa.p1.2}) follows from (\ref{prop.kappa.p1.1}) since $q \geq 2$ gives $1 \leq \log(q)/\log(2)$.
\item (\ref{prop.kappa.p1}) follows from (\ref{prop.kappa.p1.2}) since $18 \cdot 3! \cdot 2^3 = 27 \cdot 2^5$, $32^4 = 2^{20}$ and $\log(4) = 2\log(2)$.
\end{itemize}
We compute again:
\begin{align} \left\lvert \frac{\log(2)}{\log(3)} - \frac{p}{q} \right\rvert & = \frac{\lvert q\log(2) - p\log(3) \rvert}{q\log(3)} \label{prop.kappa.p2.1}
\\ & \geq \frac{q^{-27 \cdot 2^{26}\log(3)}}{q^2} \label{prop.kappa.p2.2}
\\ & \geq q^{-2^{31}} \label{prop.kappa.p2} \end{align}
This computation may be justified as follows.
\begin{itemize}
\item The equality in (\ref{prop.kappa.p2.1}) holds by rearranging.
\item (\ref{prop.kappa.p2.2}) follows from (\ref{prop.kappa.p2.1}) by exponentiating (\ref{prop.kappa.p1}) and using $\log(3) \leq q$.
\item (\ref{prop.kappa.p2}) follows from (\ref{prop.kappa.p2.2}) since $e^{9/8} \geq 3$ gives $\log(3) \leq 9/8$, hence $27 \cdot 2^{26}\log(3) + 2 \leq 27 \cdot 2^{23} \cdot 9 + 2 \leq 2^{31}$.
\end{itemize}
The inequality between the outer terms of the computation ending in (\ref{prop.kappa.p2}) verifies (\ref{prop.kappa.conc}).\\

\underline{The case $|p| \geq q + 1$}:\\

If $p \leq -1$ then $\lvert \log(2)/\log(3) - p/q \rvert \geq \log(2)/\log(3) \geq 1/2$, since $2\log(2) \geq \log(3)$. If $p \geq q + 1$ then $\lvert \log(2)/\log(3) - p/q \rvert \geq 1 - \log(2)/\log(3) \geq 1/4$, since $4\log(2) \leq 3\log(3)$. In either case (\ref{prop.kappa.conc}) holds since $q^{-2^{31}} \leq 1/4$. \\

This completes the proof of Lemma \ref{prop.kappa}. \end{Proof}

\subsection{Fixed-resolution endomorphic rigidity for $\mathbb{S}_{2,3}$ from companion article}

Define the $\times 2 \times 3$-semigroup $\mathbb{S}_{2,3} \coloneqq \{2^j3^k : j,k \in \enn\}$. The following is the specialization to $(a,b) = (2,3)$ of Clause (3) of Theorem 1.3 (Rigidity after a canonical endomorph) in \cite{RJnodyn}. The constants $\mathtt{c}_{2,3}$ and $\mathtt{s}_{2,3}$ of that clause may be evaluated as in the proof given there: we may take $\mathtt{c}_{2,3} \coloneqq 1/(8\kappa)$ and $\mathtt{s}_{2,3} \coloneqq \lceil e^{8\kappa} \rceil$ where $\kappa \in (0,\infty)_\arr$ is any exponent such that $\lvert \log(2)/\log(3) - p/q \rvert \geq q^{-\kappa}$ for all $p \in \zee$ and all $q \in [2,\infty)_\zee$.\\

Clause (3) of Theorem 1.3  in \cite{RJnodyn} carries the factor $1 - \sigma + \epsilon$ in the definition of $\ell$ and takes $t$ as the maximum of $\mathtt{s}_{2,3}$ and $\lceil e^{\ell/\mathtt{c}_{2,3}} \rceil$; we have dropped the factor and the maximum, which is permissible by the following observations.
\begin{itemize} \item The threshold on $N$ increases with $\ell$ and with $t$. 
\item The factor is at most $1$ when $\sigma \geq \epsilon$ and the conclusion is trivial when $\sigma \leq \epsilon$.
\item We have $\mathtt{s}_{2,3} \leq e^{8\kappa\ell}$ for $\ell \geq 1$.\end{itemize}

\begin{oTheorem}[Rigidity after an endomorph from $\mathbb{S}_{2,3}$] \label{othm.rigid} Let $\epsilon,\sigma \in (0,1]_\arr$, let $r \in \enn$ and let $\mu \in \prob(\rzee)$, and adopt the notation below.
\begin{align} \label{othm.rigid.h1} \ell & \coloneqq \left\lceil \frac{128}{\epsilon}\log_2\left(\frac{8(r+2)}{\epsilon}\right)\right\rceil & t & \coloneqq \bigl\lceil e^{2^{34}\ell} \bigr\rceil \end{align}
Then for all $N \in \bigl[2^{9(t+\ell)/\epsilon},\infty\bigr)_\zee$ such that $\achn(\mu,N) \geq \sigma$ there exists $m \in \mathbb{S}_{2,3} \cap [1,2N]_\zee$ such that $\einf\bigl(m_\circledcirc\mu,[0,r]_\zee\bigr) \geq \sigma - \epsilon$. \end{oTheorem}

\begin{Remark}[Invariance removes the endomorph] \label{rem.rigidinv} Since $\mathtt{T}_{ab} = \mathtt{T}_a \circ \mathtt{T}_b$, Definition \ref{def.endomorph} gives $(ab)_\circledcirc\zeta = a_\circledcirc(b_\circledcirc\zeta)$ for all $a,b \in \ennplus$ and $\zeta \in C(\rzee)^\dagger$. Thus if $\mu \in \prob(\rzee)$ is $\times 2 \times 3$-invariant, iterating Clause \ref{rem.endochangeofvar.c2} of Remark \ref{rem.endochangeofvar} gives $m_\circledcirc\mu = \mu$ for all $m \in \mathbb{S}_{2,3}$, and the conclusion of Theorem \ref{othm.rigid} becomes $\einf(\mu,[0,r]_\zee) \geq \sigma - \epsilon$. \end{Remark}

\subsection{Fourier coefficients bounded by the exponential infimum} \label{sec.rigid.master}

Lemma \ref{cor.master} is Theorem \ref{othm.rigid} read through the exponential infimum: Remark \ref{rem.foueinf} converts a lower bound on $\einf$ into an upper bound on every low-frequency Fourier coefficient, and the proof of the lemma checks that its threshold implies the one of the theorem.

\begin{Lemma}[Fourier coefficients bounded by entropy at finite resolution] \label{cor.master} For all $\epsilon \in (0,1]_\arr$, all $r \in \ennplus$, all $\times 2 \times 3$-invariant $\mu \in \prob(\rzee)$ and all $n \in \ennplus$ satisfying:
\begin{align} \label{cor.master.h1} \log\log(n) \geq \frac{2^{42}}{\epsilon} \log_2\left( \frac{8(r+2)}{\epsilon} \right) \end{align}
we have the inequality below.
\begin{align} \label{cor.master.conc} \max_{1 \leq h \leq r} \lvert \hat{\mu}(h) \rvert \leq 1 - \achn(\mu,n) + \epsilon \end{align} \end{Lemma}

The proof uses Theorem \ref{othm.rigid}, Remark \ref{rem.rigidinv} and Remark \ref{rem.foueinf}.

\begin{Proof}[of Lemma \ref{cor.master}] Let $\epsilon, r, \mu$ and $n$ be as in the statement of Lemma \ref{cor.master}, let $\ell$ and $t$ be as in (\ref{othm.rigid.h1}), and write $Z \coloneqq \log_2(8(r+2)/\epsilon)$, so that $Z \geq \log_2(24) \geq 4$. We compute:
\begin{align} \log\left(\frac{9(t+\ell)}{\epsilon}\right) & \leq \log\left(\frac{18\,e^{2^{34}\ell}}{\epsilon}\right) \label{cor.master.p1.1}
\\ & \leq 2^{34}\ell + 3 + \frac{1}{\epsilon} \label{cor.master.p1.2}
\\ & \leq \frac{2^{34} \cdot 129\,Z + 4}{\epsilon} \label{cor.master.p1.3}
\\ & \leq \frac{2^{42}}{\epsilon}\,Z \label{cor.master.p1.4}
\\ & \leq \log\log(n) \label{cor.master.p1} \end{align}
This computation may be justified as follows.
\begin{itemize}
\item (\ref{cor.master.p1.1}) holds since (\ref{othm.rigid.h1}) gives $t \leq e^{2^{34}\ell} + 1$, while $1 + \ell \leq e^{2^{34}\ell}$.
\item (\ref{cor.master.p1.2}) follows from (\ref{cor.master.p1.1}) since $\log(18) \leq 3$ and $\log(1/\epsilon) \leq 1/\epsilon$.
\item (\ref{cor.master.p1.3}) follows from (\ref{cor.master.p1.2}) since (\ref{othm.rigid.h1}) and $\epsilon \leq 1 \leq Z$ give $\ell \leq 128Z/\epsilon + 1 \leq 129Z/\epsilon$.
\item (\ref{cor.master.p1.4}) follows from (\ref{cor.master.p1.3}) since $Z \geq 4$ gives $2^{34} \cdot 129 Z + 4 \leq 2^{34} \cdot 130 Z \leq 2^{42} Z$.
\item (\ref{cor.master.p1}) follows from (\ref{cor.master.p1.4}) by (\ref{cor.master.h1}).
\end{itemize}
Exponentiating (\ref{cor.master.p1}) and using $\log(n) \leq \log_2(n)$, we obtain:
\begin{align} \label{cor.master.p2} n \geq 2^{9(t+\ell)/\epsilon} \end{align}
Now write $\sigma \coloneqq \achn(\mu,n)$. If $\sigma \leq \epsilon$ then (\ref{cor.master.conc}) holds since $|\hat{\mu}(h)| \leq 1$ for all $h \in \zee$, so we may assume $\sigma \in (\epsilon,1]_\arr$.\\

We instantiate Theorem \ref{othm.rigid} with the assignment $(\epsilon,\sigma,r,\mu,N) \coloneqq (\epsilon,\sigma,r,\mu,n)$. Under this assignment the hypothesis $N \in [2^{9(t+\ell)/\epsilon},\infty)_\zee$ agrees with (\ref{cor.master.p2}) and the hypothesis $\achn(\mu,N) \geq \sigma$ holds with equality. Thus Theorem \ref{othm.rigid} provides $m \in \mathbb{S}_{2,3} \cap [1,2n]_\zee$ with $\einf(m_\circledcirc\mu,[0,r]_\zee) \geq \sigma - \epsilon$, and since $\mu$ is $\times 2 \times 3$-invariant, Remark \ref{rem.rigidinv} gives $\einf(\mu,[0,r]_\zee) \geq \sigma - \epsilon$. Combining this with (\ref{rem.foueinf.1}) yields (\ref{cor.master.conc}). This completes the proof of Lemma \ref{cor.master}. \end{Proof}

\subsection{Decay transfer along a sequence of scales} \label{sec.rigid.scales}

Lemma \ref{lem.scales} is the bookkeeping device of the proofs of Theorems \ref{thm.retqual} and \ref{thm.eigqual}. It converts a bound that holds beyond each of a sequence of scales, together with a comparison at the next scale, into a limit superior along all scales, and it involves no measures. The proof is self-contained.

\begin{Lemma}[Decay transfer along a sequence of scales] \label{lem.scales} Let $m_1 \in \ennplus$, let $F: [m_1,\infty)_\zee \to [0,\infty)_\arr$, let $g: [m_1,\infty)_\arr \to (0,\infty)_\arr$ be nondecreasing and let $\kappa \in (0,\infty)_\arr$. Let $n_0 \in \ennplus$, and let $(J_n)_{n \in [n_0,\infty)_\zee}$ and $(\delta_n)_{n \in [n_0,\infty)_\zee}$ be sequences in $[m_1,\infty)_\arr$ and in $(0,\infty)_\arr$ respectively, with $\lim_{n \to \infty} J_n = \infty$, such that the following two conditions hold.
\begin{align} \label{lem.scales.h1} F(m) \leq \delta_n \mbox{ for all } n \in [n_0,\infty)_\zee \mbox{ and all } m \in [J_n,\infty)_\zee && \limsup_{n \to \infty}\ \frac{\log(\delta_n)}{g(J_{n+1})} \leq -\kappa \end{align}
Then we have:
\begin{align} \label{lem.scales.conc} \limsup_{m \to \infty}\ \frac{\log(F(m))}{g(m)} \leq -\kappa \end{align} \end{Lemma}

\begin{Proof}[of Lemma \ref{lem.scales}] Let $m_1,F,g,\kappa,n_0,(J_n)_{n \in [n_0,\infty)_\zee}$ and $(\delta_n)_{n \in [n_0,\infty)_\zee}$ be as in the statement of Lemma \ref{lem.scales}, and let $\epsilon \in (0,\kappa)_\arr$. By the right condition in (\ref{lem.scales.h1}) we may choose $n_1 \in [n_0,\infty)_\zee$ such that the following holds for all $n \in [n_1,\infty)_\zee$.
\begin{align} \label{lem.scales.p1} \log(\delta_n) \leq -(\kappa - \epsilon)\,g(J_{n+1}) \end{align}
Now let $m \in [m_1,\infty)_\zee$ satisfy $m \geq J_{n_1}$, and define a subset $S$ of $\enn$ as follows.
\begin{align*} S \coloneqq \{n \in [n_1,\infty)_\zee : J_n \leq m\} \end{align*}
Then $n_1 \in S$, and $S$ is finite since $\lim_{n \to \infty} J_n = \infty$. Thus $n \coloneqq \max S$ is well-defined, and $n + 1 \notin S$ gives $m \leq J_{n+1}$. We compute:
\begin{align} \log\bigl(F(m)\bigr) & \leq \log(\delta_n) \label{lem.scales.p2.1}
\\ & \leq -(\kappa - \epsilon)\,g(J_{n+1}) \label{lem.scales.p2.2}
\\ & \leq -(\kappa - \epsilon)\,g(m) \label{lem.scales.p2} \end{align}
This computation may be justified as follows.
\begin{itemize}
\item The inequality in (\ref{lem.scales.p2.1}) holds by the left condition in (\ref{lem.scales.h1}), since $n \in S$ gives $m \in [J_n,\infty)_\zee$, and since $\log$ is nondecreasing on $[0,\infty)_\arr$ under the convention \ref{conv.loginf}.
\item (\ref{lem.scales.p2.2}) follows from (\ref{lem.scales.p2.1}) by (\ref{lem.scales.p1}), since $n \geq n_1$.
\item (\ref{lem.scales.p2}) follows from (\ref{lem.scales.p2.2}) since $g$ is nondecreasing, $m \leq J_{n+1}$ and $\kappa - \epsilon \geq 0$.
\end{itemize}
Dividing (\ref{lem.scales.p2}) by the positive quantity $g(m)$ gives $\log(F(m))/g(m) \leq -(\kappa - \epsilon)$ for all $m \in [m_1,\infty)_\zee$ with $m \geq J_{n_1}$. Thus the limit superior in (\ref{lem.scales.conc}) is at most $-(\kappa - \epsilon)$, and since $\epsilon \in (0,\kappa)_\arr$ was arbitrary it is at most $-\kappa$. This completes the proof of Lemma \ref{lem.scales}. \end{Proof}

\subsection{Entropy decay from a return defect} \label{sec.rigid.cert}

Lemma \ref{lem.engine} is Lemma \ref{cor.master} at the frequencies $6^m - 1$ of Definition \ref{def.ret}, with the precision tied to the return defect. It is the one place where the threshold arithmetic of Lemma \ref{cor.master} is carried out, and the three clauses of Theorem \ref{thm.retqual} are proved from it through Lemma \ref{lem.scales}.

\begin{Lemma}[Entropy decay from a return defect] \label{lem.engine} Let $\mu \in \prob(\rzee)$ be $\times 2 \times 3$-invariant, let $m \in \ennplus$ and let $\alpha \in [2^{-m},1]_\arr$ satisfy $\ret(\mu \verr m) \leq \alpha$. Then for all $n \in \ennplus$ satisfying $\log\log(n) \geq 2^{47}m/\alpha$ we have $\achn(\mu,n) \leq 2\alpha$. \end{Lemma}

The proof below uses Clause \ref{rem.retdist.c3} of Remark \ref{rem.retdist} and Lemma \ref{cor.master}.

\begin{Proof}[of Lemma \ref{lem.engine}] Let $\mu,m,\alpha$ and $n$ be as in the statement of Lemma \ref{lem.engine}. By Clause \ref{rem.retdist.c3} of Remark \ref{rem.retdist} with $d \coloneqq m$ we may choose $h \in [1,6^m-1]_\zee$ satisfying:
\begin{align} \label{lem.engine.p1} \lvert \hat{\mu}(h) \rvert \geq 1 - \ret(\mu \verr m) \end{align}
We compute:
\begin{align} \frac{2^{42}}{\alpha}\log_2\left(\frac{8(6^m+1)}{\alpha}\right) & = \frac{2^{42}}{\alpha}\Bigl(3 + \log_2(6^m+1) - \log_2(\alpha)\Bigr) \label{lem.engine.p2.1}
\\ & \leq \frac{2^{42}}{\alpha}\Bigl(3 + 4m - \log_2(\alpha)\Bigr) \label{lem.engine.p2.2}
\\ & \leq \frac{2^{42}}{\alpha}\bigl(3 + 4m + m\bigr) \label{lem.engine.p2.3}
\\ & \leq \frac{2^{47}m}{\alpha} \label{lem.engine.p2.4}
\\ & \leq \log\log(n) \label{lem.engine.p2} \end{align}
This computation may be justified as follows.
\begin{itemize}
\item (\ref{lem.engine.p2.2}) follows from (\ref{lem.engine.p2.1}) since $\log_2(6^m+1) \leq \log_2(2 \cdot 6^m) \leq 4m$.
\item (\ref{lem.engine.p2.3}) follows from (\ref{lem.engine.p2.2}) since our assumption that $\alpha \geq 2^{-m}$ gives $-\log_2(\alpha) \leq m$.
\item (\ref{lem.engine.p2.4}) follows from (\ref{lem.engine.p2.3}) since $3 + 5m \leq 32m$.
\item (\ref{lem.engine.p2}) follows from (\ref{lem.engine.p2.4}) by the hypothesis on $n$.
\end{itemize}
We now instantiate Lemma \ref{cor.master} with the assignment $(\epsilon,r,\mu,n) \coloneqq (\alpha,6^m-1,\mu,n)$, which is admissible since $\alpha \in (0,1]_\arr$. Under this assignment the hypothesis (\ref{cor.master.h1}) holds by the computation ending in (\ref{lem.engine.p2}), so (\ref{cor.master.conc}) gives, since $h \in [1,6^m-1]_\zee$:
\begin{align} \label{lem.engine.p3} \lvert \hat{\mu}(h) \rvert \leq 1 - \achn(\mu,n) + \alpha \end{align}
Combining (\ref{lem.engine.p1}) and (\ref{lem.engine.p3}) with the hypothesis $\ret(\mu \verr m) \leq \alpha$ gives $1 - \alpha \leq 1 - \achn(\mu,n) + \alpha$, and rearranging gives $\achn(\mu,n) \leq 2\alpha$. This completes the proof of Lemma \ref{lem.engine}. \end{Proof}

\section{Proofs of theorems from Section \ref{cha.introapp.four}} \label{cha.pfret}

\subsection{Proof of Theorem \ref{thm.entret}}

Theorem \ref{thm.entret} is Lemma \ref{cor.master} at the frequency where the return defect is attained. The proof uses Clause \ref{rem.retdist.c3} of Remark \ref{rem.retdist} and Lemma \ref{cor.master}.

\begin{Proof}[of Theorem \ref{thm.entret}] Let $\epsilon,d,\mu$ and $n$ be as in the statement of Theorem \ref{thm.entret}. By Clause \ref{rem.retdist.c3} of Remark \ref{rem.retdist} there exists $h \in [1,6^d-1]_\zee$ such that the following holds.
\begin{align*} |\hat{\mu}(h)| \geq 1 - \ret(\mu \verr d) \end{align*}
We now instantiate Lemma \ref{cor.master} with the assignment $(\epsilon,r,\mu,n) \coloneqq (\epsilon,6^d-1,\mu,n)$. Under this assignment the hypothesis (\ref{cor.master.h1}) agrees with (\ref{thm.entret.h1}), so (\ref{cor.master.conc}) gives:
\begin{align*} |\hat{\mu}(h)| \leq 1 - \achn(\mu,n) + \epsilon \end{align*} 
Combining the two previous displays completes the proof of Theorem \ref{thm.entret}. \end{Proof}

\subsection{Proof of Theorem \ref{thm.retqual}}

Each clause of Theorem \ref{thm.retqual} is proved by choosing, from the right side of its inequality, a sequence of return-defect bounds $\alpha_m$ and then feeding them to Lemma \ref{lem.engine} at the scales $\exp\exp(2^{47}m/\alpha_m)$. We transfer the decay to all scales by Lemma \ref{lem.scales}. The proof uses Lemma \ref{lem.engine} and Lemma \ref{lem.scales}.

\begin{Proof}[of Theorem \ref{thm.retqual}] Let $\mu$ be as in the statement of Theorem \ref{thm.retqual}. For $r \in [16,\infty)_\arr$ we write $L(r) \coloneqq \log\log(r)$, so that $L$ is nondecreasing and $\log(L(r)) = \log^{(3)}(r)$. We also write $m_1 \coloneqq \lceil e^{e^4} \rceil$. For $r \in [m_1,\infty)_\arr$ we have $L(r) \geq 4$, hence $\log^{(3)}(r) \geq \log(4) \geq 4/3$ and $\log^{(4)}(r) \geq \log(4/3)$, so that each of $\log^{(4)}$, $\sqrt{\log^{(3)}}$ and $\log^{(3)}$ is positive and nondecreasing on $[m_1,\infty)_\arr$.\\

Each clause is verified by choosing, from the right side of its inequality, an integer $m_0 \in \ennplus$ and a sequence $(\alpha_m)_{m \in [m_0,\infty)_\zee}$ with $\alpha_m \in [2^{-m},1]_\arr$ and $\ret(\mu \verr m) \leq \alpha_m$ for all $m \in [m_0,\infty)_\zee$, and then applying Lemma \ref{lem.scales} with the function $F \coloneqq \achn(\mu,\cdot)$ on $[m_1,\infty)_\zee$ and the assignment below, where $g$ is the denominator at the left of the inequality of the clause.
\begin{align} \label{thm.retqual.p0.1} (n_0,J_m,\delta_m) \coloneqq \Bigl(m_0,\ \exp\bigl(\exp(2^{47}m/\alpha_m)\bigr),\ 2\alpha_m\Bigr) \end{align}
We record three facts about this assignment, valid for any such sequence.
\begin{itemize}
\item Since $\alpha_m \leq 1$ gives $2^{47}m/\alpha_m \geq 2^{47}$, we have $J_m \geq \exp(\exp(2^{47})) \geq m_1$, and $J_m \to \infty$ as $m \to \infty$.
\item For $r \in [J_m,\infty)_\zee$ we have $L(r) \geq L(J_m) = 2^{47}m/\alpha_m$, so Lemma \ref{lem.engine} gives $\achn(\mu,r) \leq 2\alpha_m = \delta_m$. Thus the left condition in (\ref{lem.scales.h1}) holds.
\item Since $L(J_{m+1}) = 2^{47}(m+1)/\alpha_{m+1}$, taking logarithms gives:
\begin{align} \label{thm.retqual.p0.2} \log^{(3)}(J_{m+1}) = 47\log(2) + \log(m+1) - \log(\alpha_{m+1}) \end{align}
\end{itemize}
It therefore remains, in each clause, to verify the right condition in (\ref{lem.scales.h1}). We now separate the three clauses.\\

\underline{Verification of Clause \ref{thm.retqual.c1}}:\\

Write $A$ for the limit superior at the right of (\ref{thm.retqual.h1}). Since $\ret(\mu \verr m) \leq 2$ for all $m \in \ennplus$ by Definition \ref{def.ret}, we have $A \leq 0$, and since $\achn(\mu,r) \leq 1$ for all $r \in \ennplus$ by the note following Definition \ref{def.entropy}, the left side of (\ref{thm.retqual.h1}) is at most $0$. Thus (\ref{thm.retqual.h1}) holds if $A = 0$, and we assume without loss of generality for the remainder of the current verification that $A < 0$, where $A = -\infty$ is allowed.\\

Let $s \in (0,-A)_\arr$, so we may choose $m_0 \in [3,\infty)_\zee$ with $m_0 \geq 16s^2$ such that the following holds for all $m \in [m_0,\infty)_\zee$.
\begin{align} \label{thm.retqual.p1.1} \ret(\mu \verr m) \leq \log(m)^{-s} \end{align}
For $m \in [m_0,\infty)_\zee$ we define $\alpha_m \coloneqq \log(m)^{-s}$. Since $m \geq 3$ gives $\log(m) \geq 1$, we have $\alpha_m \leq 1$. We compute:
\begin{align} s\log\log(m) & \leq s\log(m) \label{thm.retqual.p1.2.1}
\\ & \leq 2s\sqrt{m} \label{thm.retqual.p1.2.2}
\\ & \leq m\log(2) \label{thm.retqual.p1.2} \end{align}
This computation may be justified as follows.
\begin{itemize}
\item The inequality in (\ref{thm.retqual.p1.2.1}) holds since $\log(m) \leq m$ and $\log(m) \geq 1$ give $\log\log(m) \leq \log(m)$.
\item (\ref{thm.retqual.p1.2.2}) follows from (\ref{thm.retqual.p1.2.1}) since $\log(x) \leq 2\sqrt{x}$ for all $x \in [1,\infty)_\arr$.
\item (\ref{thm.retqual.p1.2}) follows from (\ref{thm.retqual.p1.2.2}) since $m \geq 16s^2$ and $\log(2) \geq 1/2$ give $\sqrt{m}\log(2) \geq 2s$.
\end{itemize}
Exponentiating (\ref{thm.retqual.p1.2}) gives $\alpha_m \geq 2^{-m}$. We apply Lemma \ref{lem.scales} with the assignment (\ref{thm.retqual.p0.1}) and $(g,\kappa) \coloneqq (\log^{(4)},s)$, for which the left condition in (\ref{lem.scales.h1}) was verified above. For the right condition, since $\log\log(m) \to \infty$ and $\log\log(m+1)/\log(m) \to 0$ as $m \to \infty$, we may fix $m_2 \in [m_0,\infty)_\zee$ such that the following hold for all $m \in [m_2,\infty)_\zee$.
\begin{align*} \log\log(m) & \geq \frac{\log(2)}{s} & 48\log(2) + s\log\log(m+1) & \leq 2\log(m) \end{align*}
Now let $m \in [m_2,\infty)_\zee$ and compute:
\begin{align} \frac{\log(\delta_m)}{g(J_{m+1})} & = \frac{\log(2) - s\log\log(m)}{\log\Bigl(47\log(2) + \log(m+1) + s\log\log(m+1)\Bigr)} \label{thm.retqual.p1.3.1}
\\ & \leq \frac{\log(2) - s\log\log(m)}{\log\Bigl(48\log(2) + \log(m) + s\log\log(m+1)\Bigr)} \label{thm.retqual.p1.3.2}
\\ & \leq \frac{\log(2) - s\log\log(m)}{\log\bigl(3\log(m)\bigr)} \label{thm.retqual.p1.3.3}
\\ & = \frac{\log(2) - s\log\log(m)}{\log\log(m) + \log(3)} \label{thm.retqual.p1.3} \end{align}
This computation may be justified as follows.
\begin{itemize}
\item The equality in (\ref{thm.retqual.p1.3.1}) holds by (\ref{thm.retqual.p0.1}), (\ref{thm.retqual.p0.2}) and the definition of $\alpha_{m+1}$, since $g(J_{m+1}) = \log\bigl(\log^{(3)}(J_{m+1})\bigr)$.
\item (\ref{thm.retqual.p1.3.2}) follows from (\ref{thm.retqual.p1.3.1}) since the numerator is nonpositive by the left condition on $m_2$, while $\log(m+1) \leq \log(2) + \log(m)$ shows that the denominator does not decrease.
\item (\ref{thm.retqual.p1.3.3}) follows from (\ref{thm.retqual.p1.3.2}) since the numerator is nonpositive, while the right condition on $m_2$ gives $48\log(2) + \log(m) + s\log\log(m+1) \leq 3\log(m)$, so that the denominator does not decrease.
\item The equality in (\ref{thm.retqual.p1.3}) holds since $\log(3\log(m)) = \log(3) + \log\log(m)$.
\end{itemize}
Since (\ref{thm.retqual.p1.3}) tends to $-s$ as $m \to \infty$, we find that the right condition in (\ref{lem.scales.h1}) holds with $\kappa = s$. Then the conclusion (\ref{lem.scales.conc}) of Lemma \ref{lem.scales} gives that the left side of (\ref{thm.retqual.h1}) is at most $-s$. Since $s \in (0,-A)_\arr$ was arbitrary, this completes the verification of Clause \ref{thm.retqual.c1}.\\

\underline{Verification of Clause \ref{thm.retqual.c2}}:\\

Write $B$ for the limit superior at the right of (\ref{thm.retqual.h3}). As with Clause \ref{thm.retqual.c1} we assume without loss of generality for the remainder of the current verification that $B < 0$, where $B = -\infty$ is allowed.\\

Let $s \in (0,-B)_\arr$, so we may choose $m_0 \in \ennplus$ with $m_0 \geq 4s^2$ such that the following holds for all $m \in [m_0,\infty)_\zee$.
\begin{align} \label{thm.retqual.p2.1} \ret(\mu \verr m) \leq e^{-s\sqrt{\log(m)}} \end{align}
For $m \in [m_0,\infty)_\zee$ we define $\alpha_m \coloneqq e^{-s\sqrt{\log(m)}}$. Since $s\sqrt{\log(m)} \geq 0$ we have $\alpha_m \leq 1$. We compute:
\begin{align} s\sqrt{\log(m)} & \leq s\sqrt{m} \label{thm.retqual.p2.2.1}
\\ & \leq m\log(2) \label{thm.retqual.p2.2} \end{align}
This computation may be justified as follows.
\begin{itemize}
\item The inequality in (\ref{thm.retqual.p2.2.1}) holds since $\log(m) \leq m$.
\item (\ref{thm.retqual.p2.2}) follows from (\ref{thm.retqual.p2.2.1}) since $m \geq 4s^2$ and $\log(2) \geq 1/2$ give $\sqrt{m}\log(2) \geq s$.
\end{itemize}
Exponentiating (\ref{thm.retqual.p2.2}) gives $\alpha_m \geq 2^{-m}$. We apply Lemma \ref{lem.scales} with the assignment (\ref{thm.retqual.p0.1}) and $(g,\kappa) \coloneqq (\sqrt{\log^{(3)}},s)$, for which the left condition in (\ref{lem.scales.h1}) was verified above. For the right condition, since $\log(m) \to \infty$ as $m \to \infty$, we may fix $m_2 \in [m_0,\infty)_\zee$ such that $\sqrt{\log(m)} \geq \log(2)/s$ for all $m \in [m_2,\infty)_\zee$. Now let $m \in [m_2,\infty)_\zee$ and compute:
\begin{align} \frac{\log(\delta_m)}{g(J_{m+1})} & = \frac{\log(2) - s\sqrt{\log(m)}}{\sqrt{47\log(2) + \log(m+1) + s\sqrt{\log(m+1)}}} \label{thm.retqual.p2.3.1}
\\ & \leq \frac{\log(2) - s\sqrt{\log(m)}}{\sqrt{48\log(2) + \log(m) + s\sqrt{\log(2) + \log(m)}}} \label{thm.retqual.p2.3} \end{align}
This computation may be justified as follows.
\begin{itemize}
\item The equality in (\ref{thm.retqual.p2.3.1}) holds by (\ref{thm.retqual.p0.1}), (\ref{thm.retqual.p0.2}) and the definition of $\alpha_{m+1}$, since $g(J_{m+1}) = \sqrt{\log^{(3)}(J_{m+1})}$.
\item (\ref{thm.retqual.p2.3}) follows from (\ref{thm.retqual.p2.3.1}) since the numerator is nonpositive by the condition on $m_2$, while $\log(m+1) \leq \log(2) + \log(m)$ shows that the denominator does not decrease.
\end{itemize}
Since (\ref{thm.retqual.p2.3}) tends to $-s$ as $m \to \infty$, we find that the right condition in (\ref{lem.scales.h1}) holds with $\kappa = s$. Then the conclusion (\ref{lem.scales.conc}) of Lemma \ref{lem.scales} gives that the left side of (\ref{thm.retqual.h3}) is at most $-s$. Since $s \in (0,-B)_\arr$ was arbitrary, the left side of (\ref{thm.retqual.h3}) is at most $B$, which completes the verification of Clause \ref{thm.retqual.c2}.\\

\underline{Verification of Clause \ref{thm.retqual.c3}}:\\

Write $C$ for the limit superior at the right of (\ref{thm.retqual.h5}), so that the right side of (\ref{thm.retqual.h5}) is $(1-C)^{-1} - 1$. As with Clauses \ref{thm.retqual.c1} and \ref{thm.retqual.c2} we assume for the remainder of the current verification that $C < 0$, where $C = -\infty$ is allowed.\\

Let $s \in (0,-C)_\arr$, so we may choose $m_0 \in \ennplus$ with $m_0 \geq 16s^2$ such that the following holds for all $m \in [m_0,\infty)_\zee$.
\begin{align} \label{thm.retqual.p3.1} \ret(\mu \verr m) \leq m^{-s} \end{align}
For $m \in [m_0,\infty)_\zee$ we define $\alpha_m \coloneqq m^{-s}$. Since $m \geq 1$ we have $\alpha_m \leq 1$. We compute:
\begin{align} s\log(m) & \leq 2s\sqrt{m} \label{thm.retqual.p3.2.1}
\\ & \leq m\log(2) \label{thm.retqual.p3.2} \end{align}
This computation may be justified as follows.
\begin{itemize}
\item The inequality in (\ref{thm.retqual.p3.2.1}) holds since $\log(x) \leq 2\sqrt{x}$ for all $x \in [1,\infty)_\arr$.
\item (\ref{thm.retqual.p3.2}) follows from (\ref{thm.retqual.p3.2.1}) since $m \geq 16s^2$ and $\log(2) \geq 1/2$ give $\sqrt{m}\log(2) \geq 2s$.
\end{itemize}
Exponentiating (\ref{thm.retqual.p3.2}) gives $\alpha_m \geq 2^{-m}$. We apply Lemma \ref{lem.scales} with the assignment (\ref{thm.retqual.p0.1}) and $(g,\kappa) \coloneqq (\log^{(3)},s/(1+s))$, for which the left condition in (\ref{lem.scales.h1}) was verified above. For the right condition, since $\log(m) \to \infty$ as $m \to \infty$, we may fix $m_2 \in [m_0,\infty)_\zee$ such that $\log(m) \geq \log(2)/s$ for all $m \in [m_2,\infty)_\zee$. Now let $m \in [m_2,\infty)_\zee$ and compute:
\begin{align} \frac{\log(\delta_m)}{g(J_{m+1})} & = \frac{\log(2) - s\log(m)}{47\log(2) + (1+s)\log(m+1)} \label{thm.retqual.p3.3.1}
\\ & \leq \frac{\log(2) - s\log(m)}{(48+s)\log(2) + (1+s)\log(m)} \label{thm.retqual.p3.3} \end{align}
This computation may be justified as follows.
\begin{itemize}
\item The equality in (\ref{thm.retqual.p3.3.1}) holds by (\ref{thm.retqual.p0.1}), (\ref{thm.retqual.p0.2}) and the definition of $\alpha_{m+1}$, since $g(J_{m+1}) = \log^{(3)}(J_{m+1})$ and $-\log(\alpha_{m+1}) = s\log(m+1)$.
\item (\ref{thm.retqual.p3.3}) follows from (\ref{thm.retqual.p3.3.1}) since the numerator is nonpositive by the condition on $m_2$, while $\log(m+1) \leq \log(2) + \log(m)$ shows that the denominator does not decrease.
\end{itemize}
Since (\ref{thm.retqual.p3.3}) tends to $-s/(1+s)$ as $m \to \infty$, we find that the right condition in (\ref{lem.scales.h1}) holds with $\kappa = s/(1+s)$. Then the conclusion (\ref{lem.scales.conc}) of Lemma \ref{lem.scales} gives that the left side of (\ref{thm.retqual.h5}) is at most $-s/(1+s) = (1+s)^{-1} - 1$. Since $s \in (0,-C)_\arr$ was arbitrary and $(1+s)^{-1}$ decreases to $(1-C)^{-1}$ as $s$ increases to $-C$, the left side of (\ref{thm.retqual.h5}) is at most $(1-C)^{-1} - 1$. This completes the proof of Clause \ref{thm.retqual.c3}, and thus completes the proof of Theorem \ref{thm.retqual}. \end{Proof}

\section{Proofs of theorems from Section \ref{cha.introapp.dioph}} \label{cha.pfspec}

\subsection{Proof of Theorem \ref{thm.eigreturn}}

Theorem \ref{thm.eigreturn} applies Lemma \ref{cor.master} at the frequency $2^j - 1$ where Lemma \ref{lem.eigreturn} makes the coefficient large; it goes through Lemma \ref{cor.master} directly rather than through Theorem \ref{thm.entret} so as to keep the threshold in terms of $2^j$ rather than $6^j$. The proof uses Definition \ref{def.eigen}, Lemma \ref{lem.eigreturn} and Lemma \ref{cor.master}.

\begin{Proof}[of Theorem \ref{thm.eigreturn}] Let $\epsilon,t,n,j,\mu,\Psi$ and $r$ be as in the statement of Theorem \ref{thm.eigreturn}, let $(\theta_k)_{k=1}^n$ be the eigenvalue list of $\Psi$, and write $\Theta \coloneqq \{\theta_1,\ldots,\theta_n\}$ and $\mathcal{V} \coloneqq \spn(\Psi)$. We note that $\mathcal{V}$ is a finite-dimensional subspace of $\sum_{\theta \in \Theta} \mathcal{E}_\mu(\theta)$, and that we have $\lVert j\theta \rVert_\zee \leq t$ for all $\theta \in \Theta$ by hypothesis.\\

Thus Lemma \ref{lem.eigreturn} applies, and since $\mu(|\mathfrak{e}^1 - \Pi_\mathcal{V}\mathfrak{e}^1|^2) = \eir(\mu \verr \Psi)$ by Definition \ref{def.eigen}, the conclusion (\ref{lem.eigreturn.conc}) reads:
\begin{align} \label{thm.eigreturn.p1} \re\,\hat{\mu}(2^j-1) \geq 1 - 2\eir(\mu \verr \Psi) - 2\pi^2t^2 \end{align}
We now instantiate Lemma \ref{cor.master} with the assignment $(\epsilon,r,\mu,n) \coloneqq (\epsilon,2^j-1,\mu,r)$. Under this assignment the hypothesis (\ref{cor.master.h1}) agrees with (\ref{thm.eigreturn.h1}), since $8(2^j - 1 + 2) = 8(2^j+1)$, so Lemma \ref{cor.master} provides:
\begin{align} \label{thm.eigreturn.p2} \max_{1 \leq h \leq 2^j-1} \lvert \hat{\mu}(h) \rvert \leq 1 - \achn(\mu,r) + \epsilon \end{align}
We compute:
\begin{align} 1 - 2\eir(\mu \verr \Psi) - 2\pi^2t^2 & \leq \re\,\hat{\mu}(2^j-1) \label{thm.eigreturn.p3.1}
\\ & \leq \lvert \hat{\mu}(2^j-1) \rvert \label{thm.eigreturn.p3.2}
\\ & \leq \max_{1 \leq h \leq 2^j-1} \lvert \hat{\mu}(h) \rvert \label{thm.eigreturn.p3.3}
\\ & \leq 1 - \achn(\mu,r) + \epsilon \label{thm.eigreturn.p3} \end{align}
This computation may be justified as follows.
\begin{itemize}
\item The inequality in (\ref{thm.eigreturn.p3.1}) holds by (\ref{thm.eigreturn.p1}).
\item (\ref{thm.eigreturn.p3.3}) follows from (\ref{thm.eigreturn.p3.2}) since $2^j - 1 \in [1,2^j-1]_\zee$.
\item (\ref{thm.eigreturn.p3}) follows from (\ref{thm.eigreturn.p3.3}) by (\ref{thm.eigreturn.p2}).
\end{itemize}
Rearranging the inequality between the outer terms of the computation ending in (\ref{thm.eigreturn.p3}) gives $\achn(\mu,r) \leq 2\eir(\mu \verr \Psi) + 2\pi^2t^2 + \epsilon$, which completes the proof of Theorem \ref{thm.eigreturn}. \end{Proof}

\subsection{Proof of Theorem \ref{thm.eigqual}}

Each clause of Theorem \ref{thm.eigqual} is a certificate: a return time for the eigenvalues, from Lemma \ref{olem.dirichlet} or from torsion, is turned by Lemma \ref{lem.eigret} into a return-defect bound beyond a sequence of scales and then Lemma \ref{lem.scales} converts this bound into a bound on the right side in the corresponding inequality of Theorem \ref{thm.retqual}.\\

The proof uses Lemma \ref{olem.dirichlet}, Clauses \ref{lem.eigret.c1} and \ref{lem.eigret.c2} of Lemma \ref{lem.eigret}, Lemma \ref{lem.scales}, and Clauses \ref{thm.retqual.c1}, \ref{thm.retqual.c2} and \ref{thm.retqual.c3} of Theorem \ref{thm.retqual}.

\begin{Proof}[of Theorem \ref{thm.eigqual}] Let $\mu$ and $(\Psi_n)_{n \in \ennplus}$ be as in the statement of Theorem \ref{thm.eigqual}. Each clause is verified by bounding the right side of the corresponding inequality of Theorem \ref{thm.retqual}, whose left side coincides with the left side of the clause. In each case Lemma \ref{lem.eigret} provides a bound $\ret(\mu \verr m) \leq \delta_n$ for all $m \in [J_n,\infty)_\zee$, and Lemma \ref{lem.scales} is applied with $m_1 \coloneqq 3$, the function $F \coloneqq \ret(\mu \verr \cdot)$ on $[3,\infty)_\zee$, and $g$ the denominator at the right of the corresponding inequality of Theorem \ref{thm.retqual}. Since $\log(3) \geq 1$, each of $\log\log$, $\sqrt{\log}$ and $\log$ is positive and nondecreasing on $[3,\infty)_\arr$. We now separate the three clauses.\\

\underline{Verification of Clause \ref{thm.eigqual.c1}}:\\

Write $A$ for the limit superior at the right of (\ref{thm.eigqual.1}). Since $\eir(\mu \verr \Psi_n) \leq 1$ for all $n \in \ennplus$ by Definition \ref{def.eigen}, we have $A \leq 0$, and since $\achn(\mu,r) \leq 1$ for all $r \in \ennplus$ by the note following Definition \ref{def.entropy}, the left side of (\ref{thm.eigqual.1}) is at most $0$. Thus (\ref{thm.eigqual.1}) holds if $A = 0$, and we assume without loss of generality for the remainder of the current verification that $A < 0$, where $A = -\infty$ is allowed.\\

Let $s \in (0,-A)_\arr$, so we may choose $n_0 \in \ennplus$ such that the following holds for all $n \in [n_0,\infty)_\zee$.
\begin{align} \label{thm.eigqual.p1.1} \eir(\mu \verr \Psi_n) \leq n^{-s} \end{align}
Now let $n \in [n_0,\infty)_\zee$, let $(\theta_k)_{k=1}^n$ be the eigenvalue list of $\Psi_n$, and define $q_n \coloneqq 4\lceil n^{s/2} \rceil$, so that $q_n \in \ennplus$ and the following two conditions hold.
\begin{align} \label{thm.eigqual.p1.2} 4n^{s/2} \leq q_n && q_n \leq 8n^{s/2} \end{align}
Here, the left condition holds since $\lceil x \rceil \geq x$, and the right condition holds since $n^{s/2} \geq 1$ and $\lceil x \rceil \leq 2x$ for all $x \in [1,\infty)_\arr$. By Lemma \ref{olem.dirichlet} applied with $q \coloneqq q_n$, we may choose $j_n \in \ennplus$ such that:
\begin{align} \label{thm.eigqual.p1.3} j_n \leq q_n^{n} && \bigvee_{k=1}^{n} \lVert j_n\theta_k \rVert_\zee \leq \frac{1}{q_n} \end{align}
Since $q_n \geq 4$ gives $1/q_n \in [0,1/2]_\arr$, Clause \ref{lem.eigret.c1} of Lemma \ref{lem.eigret} applies with the assignment $(n,j,t,\Psi) \coloneqq (n,j_n,1/q_n,\Psi_n)$, and for all $m \in [j_n,\infty)_\zee$ we compute:
\begin{align} \ret(\mu \verr m) & \leq 2\eir(\mu \verr \Psi_n) + \frac{2\pi^2}{q_n^2} \label{thm.eigqual.p1.4.1}
\\ & \leq 2n^{-s} + \frac{2\pi^2}{16}n^{-s} \label{thm.eigqual.p1.4.2}
\\ & \leq 4n^{-s} \label{thm.eigqual.p1.4} \end{align}
This computation may be justified as follows.
\begin{itemize}
\item The inequality in (\ref{thm.eigqual.p1.4.1}) holds by (\ref{lem.eigret.conc1}) under the above assignment.
\item (\ref{thm.eigqual.p1.4.2}) follows from (\ref{thm.eigqual.p1.4.1}) by (\ref{thm.eigqual.p1.1}) and the left condition in (\ref{thm.eigqual.p1.2}).
\item (\ref{thm.eigqual.p1.4}) follows from (\ref{thm.eigqual.p1.4.2}) since $\pi^2/8 \leq 2$.
\end{itemize}
We apply Lemma \ref{lem.scales} with $(g,\kappa) \coloneqq (\log\log,s)$ and the assignment below.
\begin{align} \label{thm.eigqual.p1.5} (n_0,J_n,\delta_n) \coloneqq \bigl(n_0,\ 8^nn^{sn/2},\ 4n^{-s}\bigr) \end{align}
Then $J_n \geq 8 \geq m_1$, and $J_n \to \infty$ as $n \to \infty$. The left condition in (\ref{lem.scales.h1}) holds since the left condition in (\ref{thm.eigqual.p1.3}) and the right condition in (\ref{thm.eigqual.p1.2}) give $j_n \leq q_n^n \leq J_n$, so that (\ref{thm.eigqual.p1.4}) holds for all $m \in [J_n,\infty)_\zee$. For the right condition, since $\log(n) \to \infty$ as $n \to \infty$, we may fix $n_2 \in [n_0,\infty)_\zee$ with $n_2 \geq 2$ such that the following hold for all $n \in [n_2,\infty)_\zee$.
\begin{align*} \log(n) & \geq \frac{\log(4)}{s} & \frac{s}{2}\log(n+1) & \geq 3\log(2) \end{align*}
Now let $n \in [n_2,\infty)_\zee$ and compute:
\begin{align} \frac{\log(\delta_n)}{g(J_{n+1})} & = \frac{\log(4) - s\log(n)}{\log\Bigl((n+1)\log(8) + \frac{s}{2}(n+1)\log(n+1)\Bigr)} \label{thm.eigqual.p1.6.1}
\\ & \leq \frac{\log(4) - s\log(n)}{\log\bigl((n+1)\,s\log(n+1)\bigr)} \label{thm.eigqual.p1.6.2}
\\ & \leq \frac{\log(4) - s\log(n)}{\log(n) + \log(4s) + \log\log(n)} \label{thm.eigqual.p1.6} \end{align}
This computation may be justified as follows.
\begin{itemize}
\item The equality in (\ref{thm.eigqual.p1.6.1}) holds by (\ref{thm.eigqual.p1.5}), since $\log(J_{n+1}) = (n+1)\log(8) + \frac{s}{2}(n+1)\log(n+1)$.
\item (\ref{thm.eigqual.p1.6.2}) follows from (\ref{thm.eigqual.p1.6.1}) since the numerator is nonpositive by the left condition on $n_2$, while the right condition on $n_2$ gives $(n+1)\log(8) = 3(n+1)\log(2) \leq \frac{s}{2}(n+1)\log(n+1)$, so that the denominator does not decrease.
\item (\ref{thm.eigqual.p1.6}) follows from (\ref{thm.eigqual.p1.6.2}) since the numerator is nonpositive, while $n \geq 2$ gives $n + 1 \leq 2n$ and $\log(n+1) \leq 2\log(n)$, hence $(n+1)\,s\log(n+1) \leq 4sn\log(n)$, so that the denominator does not decrease.
\end{itemize}
Since (\ref{thm.eigqual.p1.6}) tends to $-s$ as $n \to \infty$, we find that the right condition in (\ref{lem.scales.h1}) holds with $\kappa = s$. Then the conclusion (\ref{lem.scales.conc}) of Lemma \ref{lem.scales} gives that the right side of (\ref{thm.retqual.h1}) is at most $-s$, so that Clause \ref{thm.retqual.c1} of Theorem \ref{thm.retqual} bounds the left side of (\ref{thm.eigqual.1}) by $-s$. Since $s \in (0,-A)_\arr$ was arbitrary, this completes the verification of Clause \ref{thm.eigqual.c1}.\\

\underline{Verification of Clause \ref{thm.eigqual.c2}}:\\

Write $B$ for the limit superior at the right of (\ref{thm.eigqual.2}), so that the right side of (\ref{thm.eigqual.2}) is $-\sqrt{-2B}$. As with Clause \ref{thm.eigqual.c1} we assume without loss of generality for the remainder of the current verification that $B < 0$, where $B = -\infty$ is allowed.\\

Let $s \in (0,-B)_\arr$, so we may choose $n_0 \in [2,\infty)_\zee$ such that the following holds for all $n \in [n_0,\infty)_\zee$.
\begin{align} \label{thm.eigqual.p2.1} \eir(\mu \verr \Psi_n) \leq e^{-sn} \end{align}
Now let $n \in [n_0,\infty)_\zee$, let $(\theta_k)_{k=1}^n$ be the eigenvalue list of $\Psi_n$, and define $q_n \coloneqq \lceil e^{sn/2} \rceil$. Since $e^{sn/2} > 1$ we have $q_n \in [2,\infty)_\zee$ and the following two conditions hold.
\begin{align} \label{thm.eigqual.p2.2} e^{sn/2} \leq q_n && q_n \leq 2e^{sn/2} \end{align}
Here, the left condition holds since $\lceil x \rceil \geq x$, and the right condition holds since $\lceil x \rceil \leq 2x$ for all $x \in [1,\infty)_\arr$. By Lemma \ref{olem.dirichlet} applied with $q \coloneqq q_n$, we may choose $j_n \in \ennplus$ such that:
\begin{align} \label{thm.eigqual.p2.3} j_n \leq q_n^{n} && \bigvee_{k=1}^{n} \lVert j_n\theta_k \rVert_\zee \leq \frac{1}{q_n} \end{align}
Since $q_n \geq 2$ gives $1/q_n \in [0,1/2]_\arr$, Clause \ref{lem.eigret.c1} of Lemma \ref{lem.eigret} applies with the assignment $(n,j,t,\Psi) \coloneqq (n,j_n,1/q_n,\Psi_n)$, and for all $m \in [j_n,\infty)_\zee$ we compute:
\begin{align} \ret(\mu \verr m) & \leq 2\eir(\mu \verr \Psi_n) + \frac{2\pi^2}{q_n^2} \label{thm.eigqual.p2.4.1}
\\ & \leq 2e^{-sn} + 2\pi^2e^{-sn} \label{thm.eigqual.p2.4.2}
\\ & \leq 22e^{-sn} \label{thm.eigqual.p2.4} \end{align}
This computation may be justified as follows.
\begin{itemize}
\item The inequality in (\ref{thm.eigqual.p2.4.1}) holds by (\ref{lem.eigret.conc1}) under the above assignment.
\item (\ref{thm.eigqual.p2.4.2}) follows from (\ref{thm.eigqual.p2.4.1}) by (\ref{thm.eigqual.p2.1}) and the left condition in (\ref{thm.eigqual.p2.2}).
\item (\ref{thm.eigqual.p2.4}) follows from (\ref{thm.eigqual.p2.4.2}) since $2 + 2\pi^2 \leq 22$.
\end{itemize}
We apply Lemma \ref{lem.scales} with $(g,\kappa) \coloneqq (\sqrt{\log},\sqrt{2s})$ and the assignment below.
\begin{align} \label{thm.eigqual.p2.5} (n_0,J_n,\delta_n) \coloneqq \bigl(n_0,\ 2^ne^{sn^2/2},\ 22e^{-sn}\bigr) \end{align}
Then $n \geq 2$ gives $J_n \geq 4 \geq m_1$, and $J_n \to \infty$ as $n \to \infty$. The left condition in (\ref{lem.scales.h1}) holds since the left condition in (\ref{thm.eigqual.p2.3}) and the right condition in (\ref{thm.eigqual.p2.2}) give $j_n \leq q_n^n \leq J_n$, so that (\ref{thm.eigqual.p2.4}) holds for all $m \in [J_n,\infty)_\zee$. For the right condition, let $n \in [n_0,\infty)_\zee$ and compute:
\begin{align} \frac{\log(\delta_n)}{g(J_{n+1})} & = \frac{\log(22) - sn}{\sqrt{(n+1)\log(2) + \frac{s}{2}(n+1)^2}} \label{thm.eigqual.p2.6.1}
\\ & = \frac{\log(22) - sn}{(n+1)\sqrt{s/2}\,\sqrt{1 + \frac{2\log(2)}{s(n+1)}}} \label{thm.eigqual.p2.6} \end{align}
This computation may be justified as follows.
\begin{itemize}
\item The equality in (\ref{thm.eigqual.p2.6.1}) holds by (\ref{thm.eigqual.p2.5}), since $\log(J_{n+1}) = (n+1)\log(2) + \frac{s}{2}(n+1)^2$.
\item The equality in (\ref{thm.eigqual.p2.6}) holds by factoring $\frac{s}{2}(n+1)^2$ out of the square root.
\end{itemize}
Since (\ref{thm.eigqual.p2.6}) tends to $-s/\sqrt{s/2} = -\sqrt{2s}$ as $n \to \infty$, we find that the right condition in (\ref{lem.scales.h1}) holds with $\kappa = \sqrt{2s}$. Then the conclusion (\ref{lem.scales.conc}) of Lemma \ref{lem.scales} gives that the right side of (\ref{thm.retqual.h3}) is at most $-\sqrt{2s}$, so that Clause \ref{thm.retqual.c2} of Theorem \ref{thm.retqual} bounds the left side of (\ref{thm.eigqual.2}) by $-\sqrt{2s}$. Since $s \in (0,-B)_\arr$ was arbitrary, the left side of (\ref{thm.eigqual.2}) is at most $-\sqrt{-2B}$. This completes the verification of Clause \ref{thm.eigqual.c2}.\\

\underline{Verification of Clause \ref{thm.eigqual.c3}}:\\

Write $C$ for the limit superior at the right of (\ref{thm.eigqual.3}), so that the right side of (\ref{thm.eigqual.3}) is $(1-C)^{-1} - 1$. As with Clauses \ref{thm.eigqual.c1} and \ref{thm.eigqual.c2}, now with $\tor(\mu \verr q) \leq 1$ in place of $\eir(\mu \verr \Psi_n) \leq 1$, we assume for the remainder of the current verification that $C < 0$, where $C = -\infty$ is allowed.\\

Let $s \in (0,-C)_\arr$, so we may choose $q_0 \in [3,\infty)_\zee$ such that the following holds for all $q \in [q_0,\infty)_\zee$.
\begin{align} \label{thm.eigqual.p3.1} \tor(\mu \verr q) \leq q^{-s} \end{align}
Now let $n \in [q_0,\infty)_\zee$ and let $m \in [n,\infty)_\zee$. We compute:
\begin{align} \ret(\mu \verr m) & \leq 2\tor(\mu \verr m) \label{thm.eigqual.p3.2.1}
\\ & \leq 2m^{-s} \label{thm.eigqual.p3.2.2}
\\ & \leq 2n^{-s} \label{thm.eigqual.p3.2} \end{align}
This computation may be justified as follows.
\begin{itemize}
\item The inequality in (\ref{thm.eigqual.p3.2.1}) holds by (\ref{lem.eigret.conc2}) with $q \coloneqq m$.
\item (\ref{thm.eigqual.p3.2.2}) follows from (\ref{thm.eigqual.p3.2.1}) by (\ref{thm.eigqual.p3.1}), since $m \geq q_0$.
\item (\ref{thm.eigqual.p3.2}) follows from (\ref{thm.eigqual.p3.2.2}) since $m \geq n$.
\end{itemize}
We apply Lemma \ref{lem.scales} with $(g,\kappa) \coloneqq (\log,s)$ and the assignment below.
\begin{align} \label{thm.eigqual.p3.3} (n_0,J_n,\delta_n) \coloneqq \bigl(q_0,\ n,\ 2n^{-s}\bigr) \end{align}
Then $J_n = n \geq 3 = m_1$, and $J_n \to \infty$ as $n \to \infty$. The left condition in (\ref{lem.scales.h1}) is the computation ending in (\ref{thm.eigqual.p3.2}). For the right condition, let $n \in [q_0,\infty)_\zee$ and compute:
\begin{align} \label{thm.eigqual.p3.4} \frac{\log(\delta_n)}{g(J_{n+1})} = \frac{\log(2) - s\log(n)}{\log(n+1)} \end{align}
Since (\ref{thm.eigqual.p3.4}) tends to $-s$ as $n \to \infty$, we find that the right condition in (\ref{lem.scales.h1}) holds with $\kappa = s$. Then the conclusion (\ref{lem.scales.conc}) of Lemma \ref{lem.scales} gives that the right side of (\ref{thm.retqual.h5}) is at most $-s$, so that Clause \ref{thm.retqual.c3} of Theorem \ref{thm.retqual} bounds the left side of (\ref{thm.eigqual.3}) by $(1+s)^{-1} - 1$.\\

Since $s \in (0,-C)_\arr$ was arbitrary and $(1+s)^{-1}$ decreases to $(1-C)^{-1}$ as $s$ increases to $-C$, the left side of (\ref{thm.eigqual.3}) is at most $(1-C)^{-1} - 1$. This completes the proof of Clause \ref{thm.eigqual.c3}, and thus completes the proof of Theorem \ref{thm.eigqual}. \end{Proof}

\subsection{Proof of Theorem \ref{thm.atom}}

Theorem \ref{thm.atom} combines the box-average bound of Lemma \ref{lem.boxavg} with the Fourier bound of Lemma \ref{cor.master} on the same box. The proof uses Lemma \ref{lem.boxavg} and Lemma \ref{cor.master}.

\begin{Proof}[of Theorem \ref{thm.atom}] Let $\epsilon,R,\mu,\psi$ and $r$ be as in the statement of Theorem \ref{thm.atom}. We instantiate Lemma \ref{cor.master} with the assignment $(\epsilon,r,\mu,n) \coloneqq (\epsilon,2^R-1,\mu,r)$. Under this assignment, the hypothesis (\ref{cor.master.h1}) agrees with (\ref{thm.atom.h1}), so Lemma \ref{cor.master} provides:
\begin{align} \label{thm.atom.p1} \bigvee_{h=1}^{2^R-1} \lvert \hat{\mu}(h) \rvert \leq 1 - \achn(\mu,r) + \epsilon \end{align}
We compute:
\begin{align} \lvert \mu(\mathfrak{e}^1\ov{\psi}) \rvert^2 & = \lvert \langle \mathfrak{e}^1,\psi \rangle \rvert^2 \label{thm.atom.p2.1}
\\ & \leq \frac{1}{R+1} + \bigvee_{h=1}^{2^R-1} \lvert \hat{\mu}(h) \rvert \label{thm.atom.p2.2}
\\ & \leq \frac{1}{R+1} + 1 - \achn(\mu,r) + \epsilon \label{thm.atom.p2} \end{align}
This computation may be justified as follows.
\begin{itemize}
\item The equality in (\ref{thm.atom.p2.1}) holds by the form of the inner product of $L^2(\rzee,\mu)$ recorded in \ref{conv.lp}.
\item (\ref{thm.atom.p2.2}) follows from (\ref{thm.atom.p2.1}) by Lemma \ref{lem.boxavg}, which applies since $\mu$ is $\times 2$-invariant and $\psi$ is a normalized $\times 2$-eigenfunction.
\item (\ref{thm.atom.p2}) follows from (\ref{thm.atom.p2.2}) by (\ref{thm.atom.p1}).
\end{itemize}
The inequality between the outer terms of the computation ending in (\ref{thm.atom.p2}) is (\ref{thm.atom.conc}), which completes the proof of Theorem \ref{thm.atom}. \end{Proof}

\textbf{Peter Burton:} \\
Department of Mathematics and Statistics\\
University of Wyoming, Laramie WY \\
\url{peterburton1728@gmail.com}\\
\\
\textbf{Kate Juschenko:} \\
Department of Mathematics\\
University of Texas at Austin, Austin TX \\
\url{k.juschenko@austin.utexas.edu}


\begin{thebibliography}{99}

\bibitem{MR1234835} A.~Baker and G.~W\"ustholz, Logarithmic forms and group varieties, \emph{J. Reine Angew. Math.} \textbf{442} (1993), 19--62.

\bibitem{MR2563089} J.~Bourgain, E.~Lindenstrauss, P.~Michel and A.~Venkatesh, Some effective results for $\times a \times b$, \emph{Ergodic Theory Dynam. Systems} \textbf{29} (2009), 1705--1722.

\bibitem{BP26} P.~Burton and J.~Panangaden, Formulations of Furstenberg's $\times 2 \times 3$ conjecture in complex analysis and operator algebras, \emph{Ann. Funct. Anal.} \textbf{17} (2026), no.~3, Paper No.~40.

\bibitem{RJnodyn} P.~Burton and K.~Juschenko, Rudolph-Johnson without dynamics: endomorphic rigidity on the circle, preprint (2026).

\bibitem{MR0087708} J.~W.~S. Cassels, \emph{An Introduction to Diophantine Approximation}, Cambridge Tracts in Mathematics and Mathematical Physics, no.~45, Cambridge University Press, New York, 1957.

\bibitem{MR2239987} T.~M. Cover and J.~A. Thomas, \emph{Elements of Information Theory}, second ed., Wiley-Interscience, Hoboken, NJ, 2006.

\bibitem{MR3410920} T.~Eisner, B.~Farkas, M.~Haase and R.~Nagel, \emph{Operator Theoretic Aspects of Ergodic Theory}, Graduate Texts in Mathematics, vol.~272, Springer, Cham, 2015.

\bibitem{MR1681462} G.~B. Folland, \emph{Real Analysis: Modern Techniques and Their Applications}, second ed., John Wiley \& Sons, Inc., New York, 1999.

\bibitem{FGKMT18} K.~Ford, B.~Green, S.~Konyagin, J.~Maynard and T.~Tao, Long gaps between primes, \emph{J. Amer. Math. Soc.} \textbf{31} (2018), 65--105.

\bibitem{MR0213508} H.~Furstenberg, Disjointness in ergodic theory, minimal sets, and a problem in Diophantine approximation, \emph{Math. Systems Theory} \textbf{1} (1967), 1--49.

\bibitem{MR2445243} G.~H. Hardy and E.~M. Wright, \emph{An Introduction to the Theory of Numbers}, sixth ed., Oxford University Press, Oxford, 2008.

\bibitem{MR0625600} J.~E.~Hutchinson, Fractals and self-similarity, \emph{Indiana Univ. Math. J.} \textbf{30} (1981), 713--747.

\bibitem{MR1194793} A.~S.~A. Johnson, Measures on the circle invariant under multiplication by a nonlacunary subsemigroup of the integers, \emph{Israel J. Math.} \textbf{77} (1992), 211--240.

\bibitem{MR2827904} E.~Lindenstrauss, Equidistribution in homogeneous spaces and number theory, in \emph{Proceedings of the International Congress of Mathematicians, Volume~I}, Hindustan Book Agency, New Delhi, 2010, pp.~531--557.

\bibitem{Lit14} J.~E. Littlewood, Sur la distribution des nombres premiers, \emph{C. R. Acad. Sci. Paris} \textbf{158} (1914), 1869--1872.

\bibitem{Ran38} R.~A. Rankin, The difference between consecutive prime numbers, \emph{J. London Math. Soc.} \textbf{13} (1938), 242--247.

\bibitem{MR1157815} W.~Rudin, \emph{Functional Analysis}, second ed., International Series in Pure and Applied Mathematics, McGraw-Hill, New York, 1991.

\bibitem{MR1062766} D.~J. Rudolph, $\times 2$ and $\times 3$ invariant measures and entropy, \emph{Ergodic Theory Dynam. Systems} \textbf{10} (1990), 395--406.

\bibitem{MR0648108} P.~Walters, \emph{An Introduction to Ergodic Theory}, Graduate Texts in Mathematics, vol.~79, Springer-Verlag, New York--Berlin, 1982.

\end{thebibliography}
\end{document}